\documentclass[a4paper,10pt]{amsart}
\usepackage{epsfig}
\usepackage{graphicx}
\graphicspath{{~/eps}}
\usepackage[french,english]{babel}

\belowdisplayskip \abovedisplayskip
\newtheorem{theorem}{\textsc{Theorem}}[section]
\newtheorem{lemma}{\textsc{Lemma}}[section]
\newtheorem{remarque}{\textsc{Remark}}[section]

\newtheorem{definition}{\textsc{Definition}}

\newcommand{\R} {{{\rm I \! R}}}
\newcommand{\N} {{{\rm I \! N}}} 

\newcommand{\SN}{\sum_{0\leq i\leq N}}
\newcommand{\SMN}{\sum_{i,j}}
\newcommand{\SM}{\sum_{\stackrel{i=\overline{0,M}}{j=\overline{0,N}}}}
\newcommand{\sMN}{{\sum_{i,j}}^{1}}

\newcommand{\II}{{i+\frac{1}{2}}}
\newcommand{\Ii}{{i-\frac{1}{2}}}

\newcommand{\JJ}{{j+\frac{1}{2}}}
\newcommand{\Jj}{{j-\frac{1}{2}}}

\newcommand{\iJJ}{i,j+\frac{1}{2}}
\newcommand{\IIj}{i+\frac{1}{2},j}

\newcommand{\iJj}{i,j-\frac{1}{2}}
\newcommand{\Iij}{i-\frac{1}{2},j}
\newcommand{\HI}{h_{i+\frac{1}{2}}}
\newcommand{\Hi}{h_{i-\frac{1}{2}}}

\newcommand{\XII}{x_{i+\frac{1}{2}}}
\newcommand{\xii}{x_{i-\frac{1}{2}}}

\newcommand{\KJ}{k_{j+\frac{1}{2}}}

\newcommand{\lf}{\left(}
\newcommand{\rf}{\right)}

\newcommand{\Xx}{\mathcal{X}}

\newcommand{\Dh}{\Delta^{\tau}}
\newcommand{\Dhh}{\Delta^{*,\tau}}

\newcommand{\NNSU}{\Vert u \Vert_{4,\infty, \bar \Omega}}

\def\dsp{\displaystyle}

\title[Some Improvements in Finite Volume Methods ]{Some
Improvements of the Order of the Convergence of Finite Volume Solutions}
\author[B. Atfeh and A. Bradji  ]{Bilal Atfeh (1) and Abdallah BRADJI (2)*}

\begin{document}
\selectlanguage{english}

\maketitle
\begin{center}
 L.A.T.P, Universit\'e de Provence, 39 rue F.Joliot
Curie 13453 Marseille cedex 13, France\\
(1) e-mail: atfeh@cmi.univ-mrs.fr \\
(2) e-mail: bradji@cmi.univ-mrs.fr\\
*Corresponding author

\end{center}
\begin{abstract}
In this article, we improve the order of the convergence of some finite volume 
solutions approximating some second order elliptic problems.\\
In one dimensional space, we prove that finite volume approximations of order 
$O(h^{k+1})$, with $k$ integer, can be obtained after $k$ correction using the 
same scheme of three points and changing only the second members of the 
original system.\\
This is done for general smooth second order elliptic problems. These results
can be extended for non linear second order equation $ u''=f(x,u,u')$ where 
$ f$ is a smooth function .\\
In two dimensional space, we prove that finite volume approximation of order 
$O(h^2)$ can be obtained, starting with finite volume solution of order $O(h)$,
by using the same matrix and changing only the second member of the original 
system.\\
This is done for second order elliptic problems of the form $-\Delta u+pu=f$, 
with Dirichlet condition.\\
These results can be extended to obtain finite volume approximation of order 
$O(h^{k+1})$.\\  
Heart idea behind these results is the one of Fox's difference correction in 
the context of finite difference  method.
\end{abstract}
\selectlanguage{english}
\newtheorem{theo}{Theorem} 
\newtheorem{rem}{Remark}[section]
\newcommand{\noi}{\noindent}
\bibliographystyle{plain}
\selectlanguage{english}
\bibliographystyle{plain}
\section{Introduction} 
Numerical methods for partial differential equations can be divided into three 
general categories: finite difference methods, finite element methods, finite 
volume methods.\\
Finite difference and finite element methods have been attracted much more 
attention than finite volume methods, consequently there is a well developed 
literature in finite difference/finite element methods which treats several 
methods for improving the order of the convergence of the approximate solutions
those using lower scheme.\\
The desire to use low order scheme to produce highly accurate approximation in 
finite difference methods led Fox \cite{Fox} to introduce his difference 
correction technique. His idea has been modified by Pereyra and Lindberg's 
deferred correction. Theirs ideas have been developed by many authors like 
Zadunaisky's global and Frank's local defect correction (for more informations 
see \cite{Bohmer} and \cite{Skeel}).\\
Almost, the theoretical justifications of these methods are based on  the
existence of a smooth asymptotic error expansion for the base scheme. The
uniformity of the mesh and the contraction property have been the main tools to
prove such existence of the error expansions.\\
In finite Element methods, defect correction technique has been used to produce
highly order of convergence by using linear /bilinear finite element method.
This has been introduced by Barrett et al. \cite{Barrett} and Moore 
\cite{Moore} in one dimensional space by using uniform mesh and recently by 
using the so-called supraconvergent mesh condition in \cite{Butcher}.\\
In two dimensional space, under the uniform mesh Chibi \cite{Chibi} ( see also the idea of contraction property in this context in  Gao et al. \cite{Gao} ) has proved 
that, we can do only one correction on the rectangle and corrections we wish 
for periodic problems ( for a theoretical framework, you can also see the 
communication of Hackbusch in [1], pages 89-113).\\
In finite volume methods, the desire to improve the order of the convergence 
using low order scheme has not attracted the attention it merits ( see the 
introduction of \cite{Cai}). In this context, we can mention the work of Martin
et al. \cite{Martin}, where they used defect correction method, and under 
uniform mesh to propose an implicit scheme that is second order accurate both 
in time and space and uses only first jacobian for some unsteady problems.\\
The aim of this article is to develop some techniques allowing us to improve 
the order of the convergence of the finite volume solutions on arbitrary mesh 
conditions for second order elliptic problems in one and two dimensional 
spaces.\\
We prove that, starting with a finite volume solution $u^h$ of order $O(h)$ in 
$H^1$- norm, we can obtain finite volume solution of order $O(h^2)$ in 
$H^1_0$-norm, by using the same matrix that used to compute the solution $u^h$.
\\
The heart idea used in this article is the fameous Fox's difference correction 
in the context of finite difference methods.\\
The order of the convergence of the finite volume solutions, on lower schemes, 
depends on the second derivatives of the unknown solution $u$ .\\
In one dimensional space, the second derivative of $u$ can be expanded as a 
combination of the solution itself, its first derivative and a given data. We 
use this idea to obtain an optimal approximation to the second derivative by 
using the values of the basic finite volume solution $u^h$. This approximation 
allowing us to correct $u^h$ and to obtain a new approximation can be computed 
by the same matrix that used to compute $u^h$, called first correction, of 
order $O(h^2)$. \\
Other variant to compute an optimal approximation to the second derivative of 
$u$ is to use the fact that is satisfying the same equation that is satisfying 
by the solution itself but for different second member and boundary conditions 
( this holds for some second order elliptic problems). This allowing to obtain 
an optimal approximation to the second derivative, by using always the same 
matrix that used to compute $u^h$.\\
We can repeat this process, successively, to obtain finite volume approximation
s of orders $O(h^{k+1})$, where $k$ is integer by using the same matrix of 
original system.  \\
In two dimensional space, we use the second variant that used in one 
dimensional space. For the Laplacian model, the second derivatives of the 
unknown solution satisfy the same equation that is satisfying by the solution
itself, and by the same trick that used in one dimensional space, we can obtain
a new approximation to the unknown solution of order $O(h^2)$.\\
These resultes can be extended for some Dirichlet models $-\Delta u +pu=f$ and 
to obtain corrections of arbitrary order we wish. Some numerical tests 
justifying our theoretical results are done, too.\\       
\section{In One Dimentional case}
\subsection{Basic Results and Preliminaries}
The results of this article are presented in the context of classical functions
space. We denote by $C^m(\bar \Omega)$ ($\Omega $ in our paper is either an
interval in $\R$ or rectangle in $\R^2$ ) the space of continuous functions which together with
their derivatives up to order $m $ inclusive are in $C(\bar \Omega)$. The norm
is
$$
\|w\|_{m,\infty, \bar \Omega}=\max_{|\alpha|\leq
m}(\max_{\bar \Omega}|D^\alpha w|)
$$
In all that follows the letter $c$ stands for a generic, positive number,
different at each appearance but `constant` in that is independent of
discretsiation parameter $\tau, i, j$
\begin{remarque}
To show that the improvement order, will be presented (in one and two dimensional
spaces ), hold for an arbitrary admossible mesh, we try to bound each expansion
with respect to $\HI, h_i, ...$ (and we do so for two dimensional space ).
\end{remarque}   
Basic results given here are done in Eymard et al. \cite{Eymard}. Let $f$ be a 
given function defind on $(0,1)$ and consider the following equation
\begin{equation}
\label{pb1d}
\left\lbrace
\begin{array}{l}
-u_{xx}+\alpha u_x +\beta u=f(x), \,\,x \in I=(0,1), \\
\\
u(0)=u(1)=0,
\end {array}
\right.
\end{equation}
where $(\alpha ,\beta) \in \R ^{+} \times \R^{+}$.

Let $ \tau $ be an admissible mesh in the sens of \cite{Eymard}, i.e. given by 
family $(K_i)_{i=1,...,N}$, $ N \in \N ^{*}$, such that 
$K_i=(\xii  ,\XII   )$ and a family $(x_i)_{i=0,...,N+1}$
such that 
$$x_0=x_{\frac{1}{2}}=0< x_1<x_{\frac{3}{2}}<...<\xii  <x_i <
\XII  <...<x_N<x_{N+\frac{1}{2}}=x_{N+1}=1,$$
and 
$$h_i=\mbox{m }(K_i)=\XII  -\xii  , \mbox{ for }i\in
\{1,...,N\},$$
$$
h^{-}_i=x_i -\xii  \,, h^{+}_i=\XII  -x_i, 
\mbox{ for }i\in\{1,...,N\},\, h^{+}_0=h^{-}_{N+1}=0,
$$
$$\HI =x_{i+1}-x_i, i=0,...,N,\mbox{ and size}(\tau)=h=\max\{h_i, i=1,...,N\}.
$$
The system to be solved for the three points scheme by finite volume method is
\begin{equation}
\label{syslin}
\left\lbrace
\begin{array}{l}
Au^h = b, \\
u_0 = u_{N + 1} = 0,
\end{array} 
\right.
\end{equation}
where $u^h=(u_1,...,u_N)^t, b=(b_1,...,b_N)^t$ with $A $ and $ b $ are defined 
by 
\begin{equation}
\label{Mat_A}
(Au^h)_i=
\frac{1}{h_i}\lf -\frac{u_{i+1}-u_i}{\HI}+\frac{u_i-u_{i-1}}{\Hi}+\alpha(u_i-u_{i-1})\rf +\beta  u_i,\,i=1,...,N,
\end{equation}
\begin{equation}
\label{Mat_b}
b_i= \frac{1}{h_i}\int_{K_i} f(x)\,dx,\, i=1,...,N.
\end{equation}
The following theorem (see \cite{Eymard} ) gives the order of the 
convergence of the finite volume solution of the scheme of three points 
(\ref{syslin}).
\begin{theorem}
\label{TH1}
Let $f \in C([0,1]) \mbox{ and let } u \in C^2([0,1])$ be the unique solution 
of \ref{pb1d}. Let $\tau$ be an admissible mesh. Then there exists a unique
solution $u^h $ of  \ref{syslin} and the error is bounded by
\begin{equation}
\label{T1E1}
\lf \SN  \frac{(e_{i+1}-e_i)^2}{\HI }\rf^\frac{1}{2} \leq ch 
\Vert u_{2x} \Vert_{\infty,\bar I},
\end{equation}
\begin{equation}
\label{T1E2}
\vert e_i \vert \leq ch \Vert u_{2x} \Vert_{\infty}, \,\forall i\in\{1,...,N\},
\end{equation}
where $e_0=e_{N+1}=0$ and $e_i=u(x_i)-u_i$, for all $i\in \{1,...,N\}$. 
\end{theorem}
\begin{remarque}
\label{RQ1}
The uniform estimation (\ref{T1E2}) yields that the order of the convergence in
$L^2$ norm is at least $O(h)$, but numerical results shows that in general the 
order is $O(h^2)$ when $\alpha=\beta=0$, this means that
\begin{equation}
\lf\sum^N_{i=1}h_i e_i^2\rf^\frac{1}{2} \leq ch^2.
\end{equation}
\end{remarque}
Before we will be able to give general formulation of an arbitrary correction, 
we present at first the first correction and after we give the second one, 
where additional tools will be used. The general formulation of corrections can
be given later, by using the ideas of first and second correction.
\subsection{The First Correction}
By integrating both sides of equation \ref{pb1d} over each finite volume $K_i$,
we get
\begin{equation}
\label{1deq8}
-u_x(\XII  )+u_x(\xii  )+\alpha u(\XII  )-\alpha u(\xii  )+\beta \int_{K_i}u\,
dx=\int_{K_i}f\,dx,
\end{equation}
Then, the order of convergence of finite volume solution depends on the :
approximation of the flux, values $ u(\XII  )$, $ u(\xii  )$ and the integral 
$\int_{K_i} u \,dx$.

We shall use this idea combined with one of Fox to improve the order of the 
convergence of the basic solution $u^h$ on the same scheme, i.e. using the same
matrix  $A$ that used to compute the basic solution  and changing only the 
r.h.s of \ref{syslin}. 

Looking now for an expansion to the error. Assuming $u \in C^3(\bar I)$ and 
using Taylor's formula, we can get 
\begin{equation}
\label{1deq12}
\begin{array}{rl}
-\frac{u(x_{i+1})-u(x_i)}{\HI }+\frac{u(x_i)-u(x_{i-1})}{\Hi }=&
-u_x(\XII  )+u_x(\xii  )-\frac{h_{i+1}^{-}-h^{+}_i}{2}u_{2x}(\XII  )
+\frac{h_i^{-}-h^{+}_{i-1}}{2}u_{2x}(\xii)  \\
&\\
&-R^1_{i+\frac{1}{2}}+R^1_{i-\frac{1}{2}},\,\, \forall i \in \{1,...,N\}.
\end{array}
\end{equation}
where the following estimate holds
\begin{equation}
\label{1deq13}
|R^1_{i+\frac{1}{2}}| \leq c{\HI }^2 |u_{3x}|_{\infty},
\,\,\forall i \in \{0,...,N\}.
\end{equation}
On the other hand
\begin{equation}
\label{1deq14}
u(\XII  )=u(x_i)+h^{+}_i u_x+S^1_i,\,\forall i \in \{0,...,N\},
\end{equation}
where 
\begin{equation}
\label{1deq16}
|S_i^1| \leq c{h_i^{+}}^2 |u_{2x}|_{\infty},
\,\forall i \in \{0,...,N\}.
\end{equation}
Also
\begin{equation}
\label{1deq17}
\int_{K_i} u(x)dx= h_i u(x_i)+ \frac{{h^{+}_i}^2}{2}u_x(x_i)-
\frac{{h^{-}_i}^2}{2}u_x(x_{i-1}) +T^1_i,\,\forall i \in \{1,...,N\},
\end{equation}
where 
\begin{equation}
\label{1deq19}
\vert T^1_i \vert \leq ch^2h_i\|u_{2x}\|_{\infty},
\forall i \in \{1,...,N\}.
\end{equation}
Substituting terms in (\ref{1deq8}) by their expansions found in the equalities
(\ref{1deq12}),(\ref{1deq14}) and (\ref{1deq17}), we can obtain, 
$\forall i \in \{1,...,N\}$ 
\begin{eqnarray}
\label{1deq20}
-\frac{u(x_{i+1})-u(x_i)}{\HI }&+&\frac{u(x_i)-u(x_{i-1})}{\Hi }+
\alpha u(x_i)-\alpha u(x_{i-1})+\beta h_i u(x_i)=\int_{K_i} fdx \nonumber \\
&-&\frac{h_{i+1}^{-}-h^{+}_i}{2}u_{2x}(\XII  ) +
\frac{h_i^{-}-h^{+}_{i-1}}{2}u_{2x}(\xii  ) \nonumber \\
&-&\alpha h^{+}_i u_x(x_i) +\alpha h^{+}_{i-1} u_x(x_{i-1}) 
-\frac{\beta}{2}\lf{h^{+}_i}^2u_x(x_i)-{h^{-}_i}^2u_x(x_{i-1})\rf \nonumber \\ 
&-&R^1_{i+\frac{1}{2}}+R^1_{i-\frac{1}{2}}-\alpha S^1_i+\alpha S^1_{i-1}-\beta T^1_i.
\end{eqnarray}
By Substituting $u_{2x}=\alpha u_x+\beta u-f$ in the equation (\ref{1deq20}) 
we get
\begin{eqnarray}
\label{1deq21}
-\frac{u(x_{i+1})-u(x_i)}{\HI }&+&\frac{u(x_i)-u(x_{i-1})}{\Hi }+
\alpha u(x_i)-\alpha u(x_{i-1})+\beta h_i u(x_i)=\int_{K_i} fdx \nonumber \\
&-&\frac{h_{i+1}^{-}-h^{+}_i}{2}\lf \alpha u_x(\XII)+\beta u(\XII)-f(\XII)\rf 
\nonumber \\
&+&\frac{h_i^{-}-h^{+}_{i-1}}{2}\lf \alpha u_x(\xii)+\beta u(\xii)-f(\xii)\rf 
\nonumber \\
&-&\alpha h^{+}_i u_x(x_i) +\alpha h^{+}_{i-1} u_x(x_{i-1}) 
-\frac{\beta}{2}\lf (h^{+}_i)^2u_x(x_i)-(h^{-}_i)^2u_x(x_{i-1})\rf 
\nonumber \\ 
&-&R^1_{i+\frac{1}{2}}+R^1_{i-\frac{1}{2}}-\alpha S^1_i+\alpha S^1_{i-1}-\beta T^1_i.
\end{eqnarray}
After having found an appropriate expansion of the error, we can correct
the basic solution by approximating values and pointwise derivatives of the 
unknown solution $u$ in this expansion, by theirs corresponding values and 
partial values (forward approximation ) of the basic solution  $u^h$. 

The new solution $u^h_1=(u_i^1)_{i=1}^N$, called correction, obtained after  
these changes, will be defined on the same scheme, i.e. using the same matrix 
$A$ that used to compute the basic solution $u^h$, i.e. 
$u_0^1 = u_{N + 1}^1 = 0 $ and $\forall i \in \{1,...,N\}$, we have
\begin{eqnarray}
\label{1deq22}
-\frac{u^1_{i+1}-u^1_i}{\HI }&+&\frac{u^1_i-u^1_{i-1}}{\Hi }+\alpha(u^1_i-u^1_{i-1})+\beta h_i u^1_i=\int_{K_i} fdx \nonumber \\
&-&\frac{h_{i+1}^{-}-h^{+}_i}{2}\lf \alpha \partial_1 u_i+\beta u_i-f(\XII)\rf
+\frac{h_i^{-}-h^{+}_{i-1}}{2}\lf\alpha \partial_1 u_{i-1}+\beta u_{i-1}-
f(\xii)\rf \nonumber \\
&-&\alpha h^{+}_i \partial_1u_i +\alpha h^{+}_{i-1} \partial_1u_{i-1}-
\frac{\beta}{2}\lf (h^{+}_i)^2\partial_1 u_i-(h^{-}_i)^2 \partial_1u_{i-1}\rf, 
\end{eqnarray}
where $\partial_1u_i=\frac{u_{i+1}-u_i}{\HI }$. 
\subsection{The Convergence Order of the First Correction}
To analyse the convergence of the first correction, we follow the same proof 
that used for proving the order of the convergence of the basic solution. 
By subtracting (\ref{1deq21}) from (\ref{1deq22}) side by side, the error
$e_i^1=u_i^1-u(x_i)$ will be satisfied
\begin{equation}
\label{1deq23}
-\frac{e^1_{i+1}-e^1_i}{\HI }+\frac{e^1_i-e^1_{i-1}}{\Hi }+\alpha (e^1_i- e^1_{i-1})+\beta h_i e^1_i=\gamma^1_i-\gamma^1_{i-1}+\delta^1_i,
\end{equation}
where
\begin{equation}
\label{1deq24}
\gamma^1_i= -\frac{h_{i+1}^{-}-h^{+}_i}{2}\lf\alpha (\partial_1 u_i
-u_x(\XII  ))+\beta (u_i-u(\XII  )\rf 
-\alpha h^{+}_i (\partial_1 u_i-u_x(x_i)) +R^1_{i+\frac{1}{2}}+\alpha S^1_i.
\end{equation}
\begin{eqnarray}
\label{1deq25}
\delta^1_i=-\frac{\beta}{2}\lf (h^{+}_i)^2(\partial_1 u_i-u_x(x_i))
-(h^{-}_i)^2 (\partial_1 u_{i-1}-u_x(x_{i-1}))\rf +\beta T^1_i.
\end{eqnarray}
Multiplying both sides of (\ref{1deq23}) by $e_i^1$ and summing from $i=1$ to 
$i=N$, we get
\begin{eqnarray}
\label{1deq27}
-\sum^N_{i=1}\frac{e^1_{i+1}-e^1_i}{\HI }e^1_i&+&
\sum^N_{i=1}\frac{e^1_i-e^1_{i-1}}{\Hi }e^1_i+\alpha \sum^N_{i=1}
(e^1_i- e^1_{i-1})e^1_i
+\beta \sum^N_{i=1}h_i {e^1_i}^2\nonumber\\
&=&\sum^N_{i=1}\gamma^1_i e^1_i-\sum^N_{i=1}\gamma^1_{i-1} e^1_i
+\sum^N_{i=1}\delta^1_ie^1_i.
\end{eqnarray}
This gives that
\begin{equation}
\label{1deq28}
\SN \frac{(e^1_{i+1}-e^1_i)^2}{\HI}\leq \lf \SN \HI {\gamma_i}^2 
\rf^\frac{1}{2} \lf\SN \frac{(e^1_{i+1}-e^1_i)^2}{\HI }\rf^\frac{1}{2}
+\vert \sum^N_{i=1}\delta^1_ie^1_i\vert .
\end{equation} 
To estimate second term in the r.h.s of (\ref{1deq28}), we handle each term in 
(\ref{1deq25}). Indeed, we have by using inequality (\ref{1deq19}) and 
reordering sum of second term of (\ref{1deq25}), we can get
\begin{eqnarray*}
\vert \sum^N_{i=1}\delta^1_ie^1_i \vert &\leq &c\lf \sum^N_{i=1}
{h^{+}_i}^2\vert \partial_1 u_i-u_x(x_i)e_i^1\vert
+\sum^N_{i=1}{h^{-}_i}^2 \vert \partial_1 u_{i-1}-u_x(x_{i-1})e_i^1\vert+
\sum^N_{i=1} \vert T^1_ie_i^1 \vert\rf  \\
&\leq& c \lf h\sum^N_{i=1}h^{+}_i\vert \partial_1 u_i-u_x(x_i)e_i^1\vert
+h\sum^N_{i=1}h^{-}_i\vert \partial_1 u_{i-1}-u_x(x_{i-1})e_i^1\vert  
+[\sum^N_{i=1}\frac{{T^1_i}^2}{h_i}]^\frac{1}{2}
[\sum^N_{i=1}h_i {e^1_i}^2]^\frac{1}{2}\rf \\
&\leq& c\{ h[\sum^N_{i=1}h_i^{+}(\partial_1 u_i-u_x(x_i))^2]^\frac{1}{2}
[\sum^N_{i=1}h_i^{+} {e^1_i}^2]^\frac{1}{2}
+h[\sum^{N}_{i=1}h^{-}_i (\partial_1 u_{i-1}-u_x(x_{i-1}))^2]^\frac{1}{2}
[\sum^N_{i=1}h_i^{-} {e^1_i}^2]^\frac{1}{2} \\
&&+[\sum^N_{i=1}\frac{{T^1_i}^2}{h_i}]^\frac{1}{2}
[\sum^N_{i=1}h_i {e^1_i}^2]^\frac{1}{2}\}.
\end{eqnarray*}
But 
\begin{equation}
\label{1deq29}
\lf \sum^N_{i=1}h_i^{+}(\partial_1 u_i-u_x(x_i))^2\rf^\frac{1}{2}
\leq \lf \SN \HI (\partial_1 u_i-u_x(x_i))^2\rf^\frac{1}{2},
\end{equation}
and
\begin{eqnarray}
\label{1deq30}
\lf\sum^{N}_{i=1}h^{-}_i (\partial_1 u_{i-1}-u_x(x_{i-1}))^2\rf^\frac{1}{2}
&=&\lf\sum^{N-1}_{i=0}h^{-}_{i+1} (\partial_1
u_i-u_x(x_i))^2\rf^\frac{1}{2} \nonumber \\
&\leq& \lf\SN \HI (\partial_1 u_i-u_x(x_i))^2\rf^\frac{1}{2}.
\end{eqnarray}
Therefore
\begin{equation}
\label{1deq31}
\vert \sum^N_{i=1}\delta^1_ie^1_i \vert\leq c\lf h\lf \SN \HI 
(\partial_1 u_i-u_x(x_i))^2\rf^\frac{1}{2}+\lf\sum^N_{i=1}\frac{{T^1_i}^2}{h_i}\rf^\frac{1}{2}\rf
\lf \sum^N_{i=1}h_i {e^1_i}^2\rf^\frac{1}{2}.
\end{equation}
Inequality  (\ref{1deq31}) combined with discrete Poincar\'e
and triangular inequalities imply the following inequality
\begin{eqnarray}
\label{1deq32}
\lf\SN \frac{(e^1_{i+1}-e^1_i)^2}{\HI }\rf^\frac{1}{2}
&\leq&\lf\SN  \HI (\frac{h_{i+1}^{-}-h^{+}_i}{2})^2
 \alpha ^2 (\partial_1 u_i -u_x(\XII  ))^2\rf^\frac{1}{2}
 \nonumber \\
 &+&\lf\SN  \HI (\frac{h_{i+1}^{-}-h^{+}_i}{2})^2 
\beta ^2 ( u_i-u(\XII  ))^2\rf^\frac{1}{2} \nonumber \\
&+&\lf\SN \alpha^2  \HI {h_i^{+}}^2 
(\partial_1 u_i -u_x(\XII  ))^2\rf^\frac{1}{2}\nonumber \\
&+&
\lf \SN  \HI (R^1_{i+\frac{1}{2}})^2 \rf^\frac{1}{2} 
+\lf \SN  \alpha^2 \HI (S^1_i)^2 \rf^\frac{1}{2}\nonumber \\
&+&c\lf h\lf \SN \HI (\partial_1 u_i-u_x(x_i))^2\rf^\frac{1}{2}
+\lf\sum^N_{i=1}\frac{{T^1_i}^2}{h_i}\rf^\frac{1}{2}\rf
\end{eqnarray}
To estimate the r.h.s of inequality (\ref{1deq32}), we need the following 
estimates
\begin{lemma}
\label{Lemma2}
Let $u^h=(u_i)$ be the basic solution defined by (\ref{syslin}), the following
estimates hold
\begin{itemize}
\item[\bf 1.] $\dsp{\lf\SN  \HI (u_i-u(\XII  ))^2\rf^\frac{1}{2}\leq 
ch \| u_{2x}\|_{\infty,\bar I}.}$
\item[\bf 2.] $\dsp{\lf\SN  \HI (\partial_1u_i-u_x(\XII  ))^2\rf^\frac{1}{2}\leq ch 
\|u_{2x}\|_{\infty,\bar I}} .$
\item[\bf 3.] $\dsp{\lf \SN  \HI (\partial_1u_i-u_x(x_i))^2\rf^\frac{1}{2}\leq ch \|
u_{2x}\|_{\infty,\bar I}}$ 
\item[\bf 4.] $\dsp{\lf \SN  \HI  {R^1_{i+\frac{1}{2}}}^2\rf^\frac{1}{2}\leq ch^2 \|
u_{3x}\|_{\infty,\bar I}} $.
\item[\bf 5.] $\dsp{\lf \SN  \HI  {S^1_i}^2\rf^\frac{1}{2}\leq ch^2 \|
u_{2x}\|_{\infty,\bar I}}$ .
\item[\bf 6.] $\dsp{\lf \sum^N_{i=1} \frac{{T^1_i}^2}{h_i}\rf^\frac{1}{2}\leq ch^2 \|
u_{2x}\|_{\infty,\bar I}}$ . 
\end{itemize}





\end{lemma}
\begin{remarque}
Estimate {\bf{1.}} of lemma \ref{Lemma2} can be done through uniform estimate 
(\ref{T1E2}), where the estimate (\ref{T1E2}) holds for one dimension 
(see remark 2.7, page 18 in \cite{Eymard}).  The proof, we wish to present, 
holds also for the case of the finite volume scheme of five points in two 
dimensional space. \\
\end{remarque}
{\bf{Proof.}}
\begin{enumerate}
\item[\bf 1.] By triangular inequality , we have
\begin{eqnarray}
\label{1deq33}
 (\SN  \HI (u_i-u(\XII  ))^2 )^\frac{1}{2}
&\leq& (\SN  \HI (u_i-u(x_i))^2)^\frac{1}{2}
   + (\SN  \HI (u(x_i)-u(\XII  ))^2)^\frac{1}{2} \nonumber \\
 &\leq& \lf\SN  \HI e_i^2\rf^\frac{1}{2}+
 ch\Vert u_{2x}\Vert_{\infty,\bar I},
\end{eqnarray}
using the fact that $e_0=0$ and $ \SN  \HI =1$ combined with Cauchy-Schwarz 
inequality to get
\begin{eqnarray*}
 \lf\SN  \HI e_i^2\rf^\frac{1}{2}&\leq&
 \lf\SN  \frac{(e_{i+1}-e_i)^2}{\HI }\rf ^\frac{1}{2} 
  \\
 &\leq& ch \Vert u_{2x}\Vert_{\infty,\bar I} .
\end{eqnarray*} 
this with (\ref{1deq33}) imply the desired inequality 1 of lemma.
\item[\bf 2.] Using the same thecnique, yields
\begin{eqnarray*}
 \lf\SN  \HI (\partial_1 u_i-u_x(\XII  ))^2 \rf^\frac{1}{2}&\leq& \lf\SN  \HI (\partial_1 u_i-\frac{u(x_{i+1})-u(x_i)}{\HI })^2 \rf^\frac{1}{2}\\&+&\lf\SN  \HI (\frac{u(x_{i+1})-u(x_i)}{\HI } -u_x(\XII  ))^2 \rf^\frac{1}{2}\\
& \leq & ch \Vert u_{2x}\Vert_{\infty,\bar I} .   
\end{eqnarray*}
\item[\bf 3.] can be obtained as done for  {\bf{1.}} and {\bf{2.}} 
\item[\bf 4.] according to inequalty (\ref{1deq13}), we have   $ \vert R^1_{i+\frac{1}{2}} \vert \leq c h^2$, this implies the inequality 4 of the lemma.
\item[\bf 5.] according to (\ref{1deq16}), $ S^1_i$ is of order $O(h^2)$ in uniform 
norm, which implies 5 of the lemma.
\item[\bf 6.] according to (\ref{1deq19}), $\frac{{T^1_i}^2}{h_i} \leq c h^4 h_i$, the
inequality 6 of lemma will be obvious. $\Box$
\end{enumerate}

Coming back now to the lemma \ref{Lemma2}, since $h_{i+1}^{-}-h^{+}_i =O(h)$, 
then we have the following  $O(h)$ improvement.
\begin{theorem}
If the unknown solution $u $ of \ref{pb1d} belonging to $C^3 (\bar I)$. Then 
the error in the first correction defined by \ref{1deq22} is of order  $O(h^2)$
in the discrete $H^1_0$ norm, i.e
\begin{equation}
\label{1deq34}
\lf\SN \frac{(e^1_{i+1}-e^1_i)^2}{\HI }\rf^\frac{1}{2}\leq c h^2 \| u
\|_{3,\infty, \bar I},
\end{equation}
where $e^1_i=u^1_i-u(x_i)$ and $(u^1_i)_i$ are the components of the first 
correction defined by (\ref{1deq22}). 
\end{theorem}
\subsubsection{{\bf{Other Variant to Estimate the Second Derivative of Unknown
Solution}}}
For some cases, like the model $-u_{2x}+\beta u=f$, we have other possibility
to approximate the second derivative $u_{2x}$ of $u$. Indeed, $u_{2x}$  
satisfies the following equation 
$$
\left\lbrace
\begin{array}{l}
-v_{xx}+\beta v=f_{2x}(x), x \in I=(0,1), \\
v(0)=-f(0),\\
v(1)=-f(1).
\end {array}
\right.
$$
Then $u_{2x}$ satisfies the same equation that is satisfying by $u$, this
allowing us to get a finite volume approximation to $u_{2x}$, provided that
$u \in C^4(\bar I)$ (see theorem \ref{TH1}), by using the same scheme that used
to compute the basic solution $u^h$, more precisely, we use the same matrix, 
that used to compute $u^h$, to compute a finite volume approximation to  
$u_{2x}$ . This idea can be used also to compute higher order of corrections. 
\subsection{Second Correction} 
The situation in the second correction is different to that of the first
correction, because it is easy to pass from the derivative into its 
forward approximation by an order of convergence $O(h)$ (see lemma 
\ref{Lemma2}). 
To get the second correction of order $O(h^3)$, we have to look for 
approximations of first and second derivative of the unknown solution, of 
orders $O(h^2)$. That is why, we discribe how to overcome this difficuly.

Assuming that $u \in C^4(\bar I)$, by similar way to that one used to compute 
an expansion for the error (21), we can get
\begin{eqnarray}
\label{1deq35}
-\frac{u(x_{i+1})-u(x_i)}{\HI }&+&\frac{u(x_i)-u(x_{i-1})}{\Hi }+
\alpha u(x_i)-\alpha u(x_{i-1})+\beta h_i u(x_i)=\int_{K_i} fdx \nonumber \\
&-&\frac{h_{i+1}^{-}-h^{+}_i}{2}u_{2x}(\XII)-
\frac{{h_{i+1}^{-}}^2-h_{i+1}^{-}h^{+}_i +{h^{+}_i}^2}{6}u_{3x}(\XII)
\nonumber\\
&+&\frac{h_i^{-}-h^{+}_{i-1}}{2}u_{2x}(\xii  ) 
+\frac{{h_i^{-}}^2-h_i^{-}h^{+}_{i-1} +{h^{+}_{i-1}}^2}{6}u_{3x}(\xii)
\nonumber \\
&-&\alpha h^{+}_i u_x(x_i) -\alpha \frac{{h^{+}_i}^2}{2}u_{2x}(x_i)   
+\alpha h^{+}_{i-1} u_x(x_{i-1})+\alpha \frac{{h^{+}_{i-1}}^2}{2}u_{2x}(x_{i-1}) \nonumber \\
&-&\frac{\beta}{2}\lf {h^{+}_i}^2u_x(x_i)-{h^{-}_i}^2u_x(x_{i-1})\rf
-\frac{\beta}{6}\lf {h^{+}_i}^3u_{2x}(x_i)+{h^{-}_i}^3u_{2x}(x_{i-1})\rf
\nonumber \\
&+&\frac{\beta}{2}{h_i^{-}}^2\Hi u_{2x}(x_{i-1}) 
-R^2_{i+\frac{1}{2}}+R^2_{i-\frac{1}{2}}-\alpha S^2_i+\alpha S^2_{i-1}-\beta T^2_i,
\end{eqnarray} 
where
\begin{eqnarray}
\label{1deq36}
|R^2_{i+\frac{1}{2}}| \leq c \HI^3 \|u\|_{4,\infty,\bar I}, \, \,
\forall i \in \{0,...,N\},
\end{eqnarray}
\begin{eqnarray}
\label{1deq37}
|S_i^2|\leq c {h^{+}_i}^3\|u\|_{3,\infty,\bar I},\,\,\forall i \in \{0,...,N\},
\end{eqnarray}
\begin{eqnarray}
\label{1deq38}
|T^2_i| \leq c h^3h_i\|u\|_{3,\infty,\bar I},\,\,\forall i \in \{1,...,N\},
\end{eqnarray} 
In order to get correction of order $O(h^3)$ taking into account the
coefficients of the pointwise derivatives in the r.h.s of (\ref{1deq35}), we 
have to find approximations of order $O(h^2) $ to pointwise first and second 
derivative, and $O(h)$ to the pointwise third derivative in discrete 
$L^2$-norm.\\
Begining by the pointwise second derivative $u_{2x}(x_i)$, and looking for 
approxomation $u_{2x}^{h,2}=(u_{2x}^{h,2})_i, i \in \{0,...,N\}$, the idea 
that we want to suggest, is based on the use of Taylor's formula and values 
of the first correction. Indeed
\begin{eqnarray}
\label{1deq39}
u_{2x}(x_i)&=& \alpha u(x_i)+\beta u_x(x_i)+f(x_i) \nonumber \\
&=&\alpha u(x_i)+\beta \lf \frac{u(x_{i+1})-u(x_i)}{\HI}-
\frac{\HI }{2}u_{2x}(x_i)+r_{i+\frac{1}{2}}\rf +f(x_i),
\end{eqnarray}
where
\begin{equation}
\label{1deq40}
r_{i+\frac{1}{2}} \leq c \HI ^2 |u_{3x}|_{\infty}.
\end{equation}
Let $ \delta_i^h$ be the positive number $1+\frac{\beta}{2}\HI $, the equation
(\ref{1deq39}) becomes as
\begin{equation}
\label{1deq41}
\delta_i^h u_{2x}(x_i)=\alpha u(x_i)+\beta \lf \frac{u(x_{i+1})-u(x_i)}{\HI }
\rf+f(x_i)+\beta r_{i+\frac{1}{2}}.
\end{equation}
Because of the trivial inequality $ 1 \leq{\delta_i^h}\leq 1+\frac{\beta}{2} $,
we can suggest the following approximation $u_{2x}(x_i)$
\begin{equation}
\label{1deq42}
(u_{2x}^{h,2})_i=\frac{\alpha}{\delta_i^h}u_i^1+\frac{\beta}{\delta_i^h}
\lf \frac{u_{i+1}^1-u_i^1}
{\HI }\rf+ \frac{f(x_i)}{\delta_i^h},\,\forall i \in \{0,...,N\}.
\end{equation}
Looking, now, for an approximation $u_{3x}^{h,2}= (u_{3x}^{h,2})_i,\{0,...,N\}$
to pointwise third derivative. Because of 
$u_{3x}=\alpha \beta u+(\alpha^2+\beta )u_x-f_x$, it is useful to suggest the 
following approximation
\begin{equation}
\label{1deq43}
(u_{3x}^{h,2})_i=\alpha \beta u_i^1+(\alpha^2+\beta )\lf 
\frac{u_{i+1}^1-u_i^1}{\HI }\rf-f_x(x_i),\forall i \in \{0,...,N\}.
\end{equation}
\begin{remarque} We can use the approximation of the second derivative, that 
used in the first correction, to compute its approximation to obtain second
correction.
\end{remarque}
We shall prove now the following lemma
\begin{lemma}
\label{Lemma3}
If the solution $u$ of the equation (\ref{pb1d}) belonging to $C^3(\bar I)$ and
$u_{2x}^{h,2}$,$u_{3x}^{h,2}$ be the discrete expansions defined respectively
by (\ref{1deq42}) and (\ref{1deq43}). Then the following estimates hold \\
\begin{itemize}
\item[\bf 1.] $\dsp{\lf\sum_{i=0}^N \HI ((u_{2x}^{h,2})_i-u_{2x}(x_i))^2\rf^\frac{1}{2}\leq ch^2 
\Vert u \Vert_{3,\infty, \bar I}}$.
\item[\bf 2.] $\dsp{ \lf\sum_{i=0}^N \HI ((u_{3x}^{h,2})_i-u_{3x}(x_i))^2 \rf^\frac{1}{2}\leq ch 
\Vert u \Vert_{3,\infty, \bar I}.}$.
\end{itemize}
\end{lemma}
{\bf{Proof.}}
Substracting (\ref{1deq42}) from (\ref{1deq41}), to get
\begin{eqnarray*}
u_{2x}(x_i)-(u_{2x}^{h,2})_i&=&\frac{\alpha}{\delta_i^h}(u(x_i)-u_i^1)+
\frac{\beta}{\delta_i^h}\lf \frac{u(x_{i+1})-u(x_i)}{\HI }
-\frac{u_{i+1}^1-u_i^1}{\HI }\rf
+\frac{\beta r_{i+\frac{1}{2}}}{\delta_i^h}. 
\end{eqnarray*}
This implies , using triangular inequality with bound uniform of $\delta_i^h$, 
that
\begin{eqnarray}
\label{1deq44}
\lf \sum_{i=0}^N \HI ((u_{2x}^{h,2})_i-u_{2x}(x_i))^2\rf^\frac{1}{2} 
&\leq& c [\lf\sum_{i=0}^N \HI (u(x_i)-u_i^1)^2\rf^\frac{1}{2}\nonumber \\
&&+\lf\sum_{i=0}^N \HI (\frac{u(x_{i+1})-u(x_i)}{\HI }-
\frac{u_{i+1}^1-u_i^1}{\HI })^2\rf^\frac{1}{2}\nonumber \\
&&+\lf\sum_{i=0}^N \HI r_{i+\frac{1}{2}}^2\rf^\frac{1}{2}]
\end{eqnarray}
Using inequalities (\ref{1deq34}) and (\ref{1deq40}) to get the desired estimation 1 of lemma \ref{Lemma3}. \\
By the same way, we can prove the second inequality. $\Box$\\
After having acheived optimal approximations for the pointwise second and third
derivative, we look now for optimal approximations for 
$(u_{2x}(\XII  ))_i$, $(u_{3x}(\XII  ))_i$,
$(u_x(x_i)_i$, $ \forall i \in \{0,...,N\}$, we have 
$ u_{2x}(\XII  )=u_{2x}(x_i)+h_i^{+}u_{3x}(x_i)+s_i$, where 
$ \vert s_i \vert \leq c {h_i^{+}}^2 \Vert u_{4x} \Vert_\infty$, this allows us
to suggest the folowing approximation for $(u_{2x}(\XII  ))_i$ 
\begin{equation}
\label{1deq45}
u_{2x}^{i+\frac{1}{2},2}=(u_{2x}^{h,2})_i+h_i^{+}(u_{3x}^{h,2})_i,
\end{equation}
for pointwise third derivative, we can suggest the following approximation
\begin{equation}
\label{1deq46}
u_{3x}^{i+\frac{1}{2},2}=(u_{3x}^{h,2})_i,
\end{equation} 
for pointwise first derivative, we can use the trick that used for pointwise 
second derivative
\begin{equation}
\label{1deq47}
u(x_i)=\frac{u(x_{i+1})-u(x_i)}{\HI }-
\frac{\HI }{2}u_{2x}(x_i)+t_i,
\end{equation}  
where $\vert t_i \vert \leq c \HI ^2 \Vert u_{3x} \vert_\infty$,
and an approximation will be suggested as follows
\begin{equation}
\label{1deq48}
(u_x^{h,2})_i=\frac{u_{i+1}^1-u_i^1}{\HI }-
\frac{\HI }{2}(u_{2x}^{h,2})_i,
\end{equation} 
We would now prove the following lemma
\begin{lemma}
\label{Lemma4}
If the solution $u$ of equation $1$ belonging to $C^4(\bar I)$. Then the
approximations $ (u_{2x}^{i+\frac{1}{2},2})_{i=0}^N $,
$ (u_{3x}^{i+\frac{1}{2},2})_{i=0}^N $ and $ ((u_x^{h,2})_i)_{i=0}^N $ defined
respectively by the expansions (\ref{1deq45}), (\ref{1deq46}) and 
(\ref{1deq48}) satisfying the following estimates
\begin{itemize}
\item[\bf 1.] $\dsp{\lf \sum_{i=0}^N \HI (u_{2x}^{i+\frac{1}{2},2}-u_{2x}(\XII  ))^2\rf^\frac{1}{2}\leq ch^2 \Vert u \Vert_{4,\infty, \bar I}}$.
\item[\bf 2.] $\dsp{\lf\sum_{i=0}^N \HI (u_{3x}^{i+\frac{1}{2},2}-u_{3x}(\XII  ))^2\rf^\frac{1}{2}\leq ch \Vert  u \Vert_{4,\infty, \bar I}}$.
\item[\bf 3.] $\dsp{\lf\sum_{i=0}^N \HI ((u_x^{h,2})_i-u_{x}(x_i))^2 \rf^\frac{1}{2}\leq ch^2 \Vert u \Vert_{3,\infty, \bar I}}$.
\end{itemize}
\end{lemma}
{\bf{Proof.}}
\begin{enumerate}
\item[\bf 1.] We proceed as done in the proof of lemma \ref{Lemma3}.
 Using triangular inequality and equality (\ref{1deq45}), to obtain
\begin{eqnarray*}
\lf \sum_{i=0}^N \HI (u_{2x}^{i+\frac{1}{2},2}-u_{2x}(\XII  ))^2\rf ^\frac{1}{2}&\leq& 
\lf\sum_{i=0}^N \HI ((u_{2x}^{h,2})_i-u_{2x}(x_i))^2\rf^\frac{1}{2} \\
&&+\lf
\sum_{i=0}^N\HI {h_i^{+}}^2((u_{3x}^{h,2})_i-u_{3x}(x_i))^2\rf^\frac{1}{2}\\
&&+\lf\sum_{i=0}^N \HI s_i^2\rf^\frac{1}{2}
\end{eqnarray*}
Using lemma \ref{Lemma3} and the uniform bound of $ s_i $ to obtain the desired 
inequality 1 of the lemma \ref{Lemma4}.
\item[\bf 2.] and {\bf{3.}} of the lemma can be handled by the same way as done for the first estimation. $ \Box $
\end{enumerate}

Now we are able to define the second correction $ u^h_2=(u_i^2)_{i=0}^{N+1} $,
where  $u_0^2=u_{N+1}^2=0 $ and for all $ i \in \{1,...,N\} $, we have
\begin{eqnarray}
\label{1deq49}
-\frac{u_{i+1}^2-u_i^2}{\HI }&+&\frac{u_i^2-u_{i-1}^2}{\Hi }+
\alpha u_i^2-\alpha u_{i-1}^2+\beta h_i u_i^2=\int_{K_i} fdx \nonumber \\
&-&\frac{h_{i+1}^{-}-h^{+}_i}{2}u_{2x}^{{i+\frac{1}{2}},2}-
\frac{{h_{i+1}^{-}}^2-h_{i+1}^{-}h^{+}_i
+{h^{+}_i}^2}{6}u_{3x}^{{i+\frac{1}{2}},2} \nonumber \\
&+&\frac{h_i^{-}-h^{+}_{i-1}}{2}u_{2x}^{{i-\frac{1}{2}},2} 
+\frac{{h_i^{-}}^2-h_i^{-}h^{+}_{i-1} +{h^{+}_{i-1}}^2}{6}
u_{3x}^{{i-\frac{1}{2}},2}\nonumber \\
&-&\alpha h^{+}_i (u_x^{h,2})_i -\alpha \frac{{h^{+}_i}^2}{2}(u_{2x}^{h,2})_i  
+\alpha h^{+}_{i-1} (u_x^{h,2})_{i-1}+\alpha \frac{{h^{+}_{i-1}}^2}{2}
(u_{2x}^{h,2})_{i-1} \nonumber \\
&-&\frac{\beta}{2}\{{h^{+}_i}^2(u_x^{h,2})_i-{h^{-}_i}^2(u_x^{h,2})_{i-1}\}
-\frac{\beta}{6}\{{h^{+}_i}^3(u_{2x}^{h,2})_i+{h^{-}_i}^3(u_{2x}^{h,2})_{i-1}\}\nonumber \\
&+&\frac{\beta}{2}{h_i^{-}}^2\Hi (u_{2x}^{h,2})_{i-1} 
\end{eqnarray} 
To analyse the error of the convergence, we proceed as done for the basic
solution and the first correction. Indeed, let $e_i^2=u_i^2-u(x_i)$ and
substracting equality (\ref{1deq35}) from (\ref{1deq49}) to get
\begin{equation}
\label{1deq50}
-\frac{e_{i+1}^2-e_i^2}{\HI }+\frac{e_i^2-e_{i-1}^2}{\Hi }+
\alpha e_i^2-\alpha e_{i-1}^2+\beta h_i e_i^2= \gamma_i^2-\gamma_{i-1}^2+\delta_i^2.
\end{equation}
where
\begin{eqnarray}
\label{1deq51}
\gamma_i^2&=&-\frac{h_{i+1}^{-}-h^{+}_i}{2}
(u_{2x}^{i+\frac{1}{2},2}-u_{2x}(\XII  ))
-\frac{{h_{i+1}^{-}}^2-h_{i+1}^{-}h^{+}_i +{h^{+}_i}^2}{6}
(u_{3x}^{{i+\frac{1}{2}},2}-u_{3x}(\XII  ))\nonumber \\
&-&\alpha h^{+}_i ((u_x^{h,2})_i-u_x(x_i))
 -\alpha \frac{{h^{+}_i}^2}{2}
((u_{2x}^{h,2})_i -u_{2x}(x_i))
 + R^2_{i+\frac{1}{2}}+\alpha S^2_i
\end{eqnarray}
 
\begin{eqnarray}
\label{1deq52}
\delta_i^2&=&-\frac{\beta}{2}\{{h^{+}_i}^2((u_x^{h,2})_i-u_x(x_i))-
{h^{-}_i}^2((u_x^{h,2})_{i-1}-u_x(x_{i-1}))\}\nonumber \\
&-&\frac{\beta}{6}\{{h^{+}_i}^3(u_{2x}^{h,2})_i-u_{2x}(x_i))+{h^{-}_i}^3
(u_{2x}^{h,2})_{i-1}-u_{2x}(x_{i-1}))\}\nonumber \\
&-&\frac{\beta}{2}{h_i^{-}}^2\Hi ((u_{2x}^{h,2})_{i-1}-u_{2x}(x_{i-1}))  
+\beta T^2_i
\end{eqnarray}
Using the proof of the convergence of the basic solution and the first
correction together with triangular and discrete Poincare inequalities 
combined with lemma \ref{Lemma4} and expansions (\ref{1deq36}), (\ref{1deq37}) and
(\ref{1deq38}) to get the following 
$O(h^2) $ improvement.
\begin{theorem}
If the unknown solution of (\ref{pb1d}) belonging to $C^4 (\bar I)$. Then the 
error in the second correction defined by (\ref{1deq49}) is of order  $O(h^3)$ 
in the discrete $H^1_0$ norm, i.e.
\begin{equation}
\label{1deq53}
\lf\SN \frac{(e^2_{i+1}-e^2_i)^2}{\HI }\rf^\frac{1}{2}\leq c h^3 \Vert u
\Vert_{4,\infty, \bar I}
\end{equation}
where $e^2_i=u^2_i-u(x_i)$ and $(u^2_i)_i$ are the components of the second 
correction defined by (\ref{1deq49}). 
\end{theorem} 
\subsection{Corrections of Higher Order}
In this section, we give the general formulation of an arbitrary correction. 
The pointwise derivatives will be approximated in the light of ones of the 
first and second correction. The proof of the order of the convergence is the 
same one that done for the first and second correction.\\
For each integer $ k \geq 1$
\begin{eqnarray*}
-\frac{u(x_{i+1})-u(x_i)}{\HI }+\frac{u(x_i)-u(x_{i-1})}
{\Hi }&=&-u_x(\XII  )+u_x(\xii )+
\sum^{k+1}_{m=2}a^i_m u_{mx}(\XII  )\\
&-&\sum^{k+1}_{m=2}a^{i-1}_m u_{mx}(\xii  )-
R^k_{i+\frac{1}{2}}+R^k_{i-\frac{1}{2}}
\end{eqnarray*}
\begin{equation}
\label{1deq55} 
a^i_m=\frac{1}{m!}\lf\sum^{m-1}_{j=0}{(h^{{-}}_{i+1})}^j{(-h^{+}_i)}^{m-1-j}\rf
.
\end{equation}
and 
\begin{equation}
\label{1deq54}
|R_\II |\leq ch_\II^{k+1}\|u\|_{k+2,\infty,\bar I}  
\end{equation}
We have also
\begin{equation}
\label{1deq56}
u(\XII  )=u(x_i)+\sum^k_{m=1} \frac{(h^{+}_i)^m}{m!} u_{mx}(x_i)+S^k_i,
\end{equation}
where 
\begin{equation}
\label{1deq57}
|S_i^k|\leq ch_i^{k+1}\|u\|_{k+1,\infty,\bar I}
\end{equation}
 For the fifth term in (\ref{1deq8}), we have
\begin{equation}
\label{1deq58}
\int_{K_i} u(x)dx= h_i u(x_i)+ 
\sum^k_{m=1} \frac{(h^{+}_i)^{m+1}}{(m+1)!}u_{mx}(x_i)
-\sum^k_{m=1}\frac{(-h^{-}_i)^{m+1}}{(m+1)!}\sum^{k-m}_{j=0}
\frac{\Hi ^j}{j!}u_{(m+j)x}(x_{i-1})
+T^k_i,
\end{equation}
where 
\begin{equation}
\label{1deq59}
|T^k_i|\leq ch_ih^{k+1}\|u\|_{k+1,\infty,\bar I}
\end{equation}
Substituting terms of (\ref{1deq8}) by theirs expansions (\ref{1deq56}),
(\ref{1deq57}) and (\ref{1deq59}), we obtain 
\begin{eqnarray}
\label{1deq60}
-\frac{u(x_{i+1})-u(x_i)}{\HI }+
\frac{u(x_i)-u(x_{i-1})}{\Hi }&+&
\alpha u(x_i)-\alpha u(x_{i-1})+\beta h_i u(x_i)=
\int_{K_i} fdx \nonumber \\
&-&\sum^{k+1}_{m=2}a^i_m u_{mx}(\XII  )+
\sum^{k+1}_{m=2}a^{i-1}_m u_{mx}(\xii  ) \nonumber\\
&-&\alpha \sum^k_{m=1} 
\frac{(h^{+}_i)^m}{m!}u_{mx}(x_i)+\alpha \sum^k_{m=1} 
\frac{(h^{+}_{i-1})^m}{m!}u_{mx}(x_{i-1})\nonumber \\
&-&\beta\{ \sum^k_{m=1} \frac{(h^{+}_i)^{m+1}}{(m+1)!}u_{mx}(x_i)
-\sum^k_{m=1}\frac{(-h^{-}_i)^{m+1}}{(m+1)!}\sum^{k-m}_{j=0}\frac{\Hi ^j}{j!}
 u_{(m+j)x}(x_{i-1}) \}\nonumber  \\
&-&R^k_{i+\frac{1}{2}}+R^k_{i-\frac{1}{2}}-
\alpha S^k_i+\alpha S^k_{i-1}-\beta T^k_i
\end{eqnarray}
After having found an expansion approximating the equation (\ref{1deq8}), we
need now the following useful lemma
\begin{lemma}
\label{Lemma5}
Each $kth$ derivative of the solution $u$ of the problem \ref{pb1d} can be 
expanded as a linear combination of the solution itself, its derivative and 
the derivatives of the given function $f$ up to and including $(k-2)th$ 
derivative, i.e, there exist  reals 
$\{\alpha^k_j\}^{k-2}_{j=0} \cup \{\overline \alpha^k_1,\overline \alpha^k_2\}$
such that
\begin{equation}
\label{1deq61}
u_{kx}=\sum ^{k-2}_{j=0} \alpha ^k_j f_{jx}+\overline \alpha^k_1 u +\overline \alpha^k_2 u_x.
\end{equation}   
\end{lemma}
{\bf{Proof. }}We can prove this lemma by induction on the integer $k$. $\Box$\\
Assuming, now, that we have obtained the $(k-1)th$ correction 
$u^h_{k-1}=(u_i^{k-1})_i$, i.e. approximation of order $O(h^k)$. According to 
equality (\ref{1deq60}), to obtain correction of order $O(h^{k+1})$,
 we have to find approximations for the pointwise derivative up and including 
$k+1$ order of the solution $u$. The idea which we will present is similar to 
that one presented  to compute second correction. 

At first, we look for optimal approximations to $ u_{2x}(x_i)$,...,
$u_{(k+1)x}(x_i)$. To do so, we use the previous correction, i.e. $(k-1)th$
correction, and the optimal approximations to $u_{2x}(x_i)$,...,$u_{kx}(x_i)$ 
used to define this correction. That is why, we define the $kth$ correction by 
induction, we assume that, we have obtained $(k-1)th$ correction of order 
$O(h^k)$ and we have found optimal approximations (according to theirs 
coefficients in (\ref{1deq60}) $(u_{2x}^{h,k-1})_i$,...,$(u_{kx}^{h,k-1})_i$ 
for $u_{2x}(x_i)$,..., $u_{kx}(x_i)$, i.e. theirs orders of convergence are 
$O(h^{k-1})$,...,$O(h)$.
 
Because the coefficient of $u_{(k+1)x}(x_i)$ in (\ref{1deq60}) is of order 
$O(h^k)$, it suffices to approximate it by in order $ O(h) $. This, can be 
done easily through lemma \ref{Lemma5},i.e. an approximation defined by 
\begin{equation}
\label{1deq62}
(u_{(k+1)x}^{h,k})_i=\sum ^{k-1}_{j=0} \alpha ^{k+1}_j f_{jx}(x_i)+
\bar \alpha^{k+1}_1 u_i^{k-1} +\bar \alpha^{k+1}_2 \partial u_i^{k-1},
\forall i \in \{0,...,N\}.
\end{equation}
 
We can use, also, in (\ref{1deq62}) instead of the the $(k-1)$ correction, 
the basic solution $u^h$.\\
For any integer $\beta$ such that $2 \leq \beta \leq k$, we look to find 
approximation $ u_{\beta x}^{h,k}$ of order $O(h^{k+2-\beta})$ to pointwise 
derivative $ (u_{\beta x}(x_i))_i$ of order $\beta$, bacause the coefficients 
of such derivative in (\ref{1deq60}) are of order $\beta -1$. We have through 
lemma \ref{Lemma5}
\begin{eqnarray}
\label{1deq63}
u_{\beta x}(x_i)&=&\sum ^{\beta-2}_{j=0} \alpha ^{\beta}_j f_{jx}(x_i)+
\bar \alpha^{\beta}_1 u(x_i) +\bar \alpha^{\beta}_2 u_x( u_i) \nonumber \\
&=&\sum ^{\beta-1}_{j=0} \alpha ^{\beta}_j f_{jx}(x_i)+\bar \alpha^{\beta}_1 u(x_i)\nonumber \\
&+&\bar \alpha^{\beta}_2 \lf \frac{u(x_{i+1})-u(x_i)}{\HI }-
\sum_{j=2}^{k-\beta+2} \frac{u_{jx}(x_i)}{j!}\HI ^{j-1}\rf+r_i^k,\,
 \forall i \in \{0,...,N\} ,
\end{eqnarray}  
where
\begin{equation}
\label{1deq64}
\vert r_i^k \vert \leq c\HI ^{k+2-\beta}\vert u \vert_{k+3-\beta,\infty}
\end{equation}    
An obvious approximation for pointwise derivative of order $\beta$ can be given
as
\begin{equation}
\label{1deq65}
(u_{\beta x}^{h,k})_i=\sum ^{\beta-2}_{j=0} \alpha ^{\beta}_j f_{jx}(x_i)
+\bar \alpha^{\beta}_1 u_i^{k-1} 
+\bar \alpha^{\beta}_2 \lf \frac{u_{i+1}^{k-1}-u_i^{k-1}}{\HI }-
\sum_{j=2}^{k-\beta+2} \frac{(u_{jx}^{h,k-1})_i}{j!}\HI ^{j-1}\rf,
\end{equation}
We would prove the following lemma
\begin{lemma}
If the solution $u$ of the equation \ref{pb1d} belonging to $C^{k+1}(\bar I )$.
Then the approximations $ u_{\beta x}^{h,k}$, where $2 \leq \beta \leq k+1$, 
defined by the expansions (\ref{1deq62}) and (\ref{1deq65}) satisfying the 
following estimate
\begin{equation}
\label{1deq66}
\lf \sum_{i=0}^N \HI ((u_{\beta x}^{h,k})_i-u_{\beta x}(x_i))^2
\rf^\frac{1}{2} \leq c h^{k-\beta+2} \Vert u \Vert_{k+1,\infty,\bar I}. 
\end{equation} 
\end{lemma}
{\bf{Proof .}} We can prove this by induction. $\Box$\\
After having found optimal approximations to fundamental pointwise derivatives,
we derive now optimal approximations for $ (u_{\beta
x}(\XII  ))_{i=0}^N$,$2 \leq \beta \leq k+1$ and $
(u_x(x_i))_{i=0}^N$.\\
We have
\begin{equation}
\label{1deq67} 
u_{\beta x}(\XII  )=\sum_{j=0}^{k-\beta+1}\frac{(h_i^{+})^j}{j!} u_{(\beta
+j)x}(x_i)+s_i^k,\forall i \in \{0,...,N\}
\end{equation} 
where
\begin{equation}
\label{1deq68} 
\vert s_i^k \vert \leq c(h_i^{+})^{k-\beta+2} \vert u_{(k+2)x} \vert_{\infty}
\end{equation} 
We can suggest the following approximation
\begin{equation}
\label{1deq69}
(u_{\beta x}^{\II,k})_i = \sum_{j=0}^{k-\beta+1}\frac{(h_{+})^j}{j!} (u_{(\beta
+j)x}^{h,k})_i,\forall i \in \{0,...,N\}.
\end{equation}
For the pointwise first derivative, we can do
\begin{equation}
\label{1deq70}
u_x(x_i)=\frac{u(x_{i+1})-u(x_i)}{\HI }-\sum_{j=2}^k 
\frac{\HI ^{j-1}}{j!} u_{jx}(x_i)+t_i^k,
\end{equation} 
where
\begin{equation}
\label{1deq71} 
\vert t_i^k \vert \leq c \HI ^k \vert u_{(k+1)x} \vert_{\infty},
\forall i \in \{0,...,N\}
\end{equation} 
this allowing us to consider the following approximation
\begin{equation}
\label{1deq72}
(u_x^{h,k})_i = \frac{u^{k-1}_{i+1}-u^{k-1}_i}{\HI }-\sum_{j=2}^k 
\frac{\HI ^{j-1}}{j!} (u_{jx}^{h,k})_i,\forall i \in \{0,...,N\}
\end{equation}
We need the following useful lemma
\begin{lemma}
\label{Lemma7}
If the solution $u$ of the equation (\ref{pb1d}) belonging to 
$C^{k+1}(\bar I )$. Then the approximations $(u_{\beta x}^{\II,k})_i$, where 
$2 \leq \beta \leq k+1$, defined and $(u_x^{h,k})_i$ defined respectively by 
the expansions (\ref{1deq69}) and (\ref{1deq72}) 
satisfying the following estimation 
$$\lf \sum_{i=0}^N \HI ((u_{\beta x}^{i+\frac{1}{2},k})_i-
u_{\beta x}(\XII  ))^2\rf^\frac{1}{2} \leq c h^{k-\beta+2} 
\Vert u \Vert_{k+1,\infty,\bar I}.$$
$$\lf\sum_{i=0}^N \HI ((u_x^{h,k})_i-u_x(x_i))^2\rf^\frac{1}{2}
 \leq c h^k \Vert u \Vert_{k+1,\infty,\bar I}.$$
\end{lemma}
{\bf{Proof. }} The proof can be done as done for proving lemma \ref{Lemma4}. $\Box$\\
Now we are able to define the $kth$ correction $u^k_h=(u_i^k)_{i=0}^{N+1}$, 
where $u_0^k=u_{N+1}^k=0$ and for $ i \in \{1,...,N\}$, we have 
\begin{eqnarray}
\label{1deq73}
-\frac{u^k_{i+1}-u^k_i}{\HI }&+&
\frac{u^k_i-u^k_{i-1}}{\Hi }+
\alpha u^k_i-\alpha u^k_{i-1}+\beta h_i u_i^k=
\int_{K_i} fdx \nonumber \\
&&-\sum^{k+1}_{m=2}a^i_m (u_{mx}^{i+\frac{1}{2},k})_i+
\sum^{k+1}_{m=2}a^{i-1}_m u_{mx}^{i-\frac{1}{2},k} 
-\alpha \sum^k_{m=1} 
\frac{(h^{+}_i)^m}{m!}(u_{mx}^{h,k})_i+\alpha \sum^k_{m=1} 
\frac{(h^{+}_{i-1})^m}{m!}(u_{mx}^{h,k})_{i-1})\nonumber \\
&&-\beta\lf \sum^k_{m=1} \frac{(h^{+}_i)^{m+1}}{(m+1)!}(u_{mx}^{h,k})_i
-\sum^k_{m=1}\frac{(-h^{-}_i)^{m+1}}{(m+1)!}\sum^{k-m}_{j=0}\frac{\Hi ^j}{j!}
 (u_{(m+j)x})^{h,k}_{i-1} \rf
\end{eqnarray}
Consedering the following expansions
\begin{equation}
\label{1eq74}
\gamma_i^k=-\sum^{k+1}_{m=2}a^i_m ((u_{mx}^{i+\frac{1}{2},k})_i-
u_{m x}(\XII  ))-\alpha \sum^k_{m=1} 
\frac{(h^{+}_i)^m}{m!}((u_{mx}^{h,k})_i-u_{mx}(x_i))
+R^k_{i+\frac{1}{2}}+\alpha S^k_i
\end{equation}

\begin{equation}
\label{1deq75}
\delta_i^k=-\beta\lf \sum^k_{m=1} \frac{(h^{+}_i)^{m+1}}{(m+1)!}
\lf(u_{mx}^{h,k})_i-u_{mx}(x_i)\rf -\sum^k_{m=1}\frac{(-h^{-}_i)^{m+1}}{(m+1)!}\sum^{k-m}_{j=0}\frac{\Hi ^j}{j!} 
\lf(u_{(m+j)x}\rf^{h,k}_{i-1}-u_{(m+j)x}(x_{i-1}) \rf+\beta T^k_i. 
\end{equation}
Let $ e_i^k=u_i^k-u(x_i)$ be the error in the $kth$ correction, thus
\begin{equation}
\label{1deq76}
-\frac{e^k_{i+1}-e^k_i}{\HI }+
\frac{e^k_i-e^k_{i-1}}{\Hi }+
\alpha e^k_i-\alpha e^k_{i-1}+\beta h_i e_i^k=\gamma_i^k-\gamma_{i-1}^k+
\delta_i^k
\end{equation}
Using the proof of the convergence of the basic solution, first and second 
corrections and lemma \ref{Lemma7} together with  (\ref{1deq55}) ,(\ref{1deq54}),
(\ref{1deq57}) and (\ref{1deq59}) to 
get the theorem
\begin{theorem}
If the unknown solution of (\ref{pb1d}) belonging to $C^{k+2} (\bar I)$. Then 
the error in the $kth$ correction, defined by (\ref{1deq76}) is, of order  
$O(h^{k+1})$ in the discrete $H^1_0$ norm, i.e.
\begin{equation}
\label{1deq77}
\lf\SN \frac{(e^k_{i+1}-e^k_i)^2}{\HI}\rf^\frac{1}{2}
\leq c h^{k+1} \vert u \vert_{k+2,\infty, \bar I}
\end{equation}
where $e^k_i=u^k_i-u(x_i)$ and $(u^k_i)_i$ are the components of the $kth$
correction defined
by (\ref{1deq73}). 
\end{theorem}
\begin{remarque}
As you have seen, we can  generalize the results obtained for general equation
$y''=f(x,y,y')$, where $f$ is a smooth function.
\end{remarque}

\newpage
\section{In Two Dimension Space}
\subsection{ Basic Results}
Considering the second order elliptic problem, with homogeneous boundary 
conditions
\begin{equation}
\label{pb2d} 
\left\lbrace
\begin{array}{l}
-\Delta u=f, \mbox{on } \Omega =(0,1)^2\\
u|_ \Gamma=0.
\end {array}
\right.
\end{equation}
where $\Gamma=\partial \Omega$ is the boundary of $\Omega$ and assuming that 
the solution $u$ is belonging to $C^2(\bar \Omega)$, and the second member 
$f \in C(\bar \Omega)$.

Let $ \tau=(K_{ij})_{1\leq i\leq M; 1\leq j\leq N}$ be an admissible mesh of
$\Omega$ in the sense of \cite{Eymard}, that is satisfying the following
assumption

{\bf{Assumption. }} {\it{Let $M,\,N \in \N ^*$, $(h_i)_{i=1}^M,\,(k_j)_{j=1}^N$ 
are positive numbers and such that
$$ \sum_{i=1}^M h_i=\sum_{j=1}^N k_j=1,$$
and let $h_0=h_{M+1}=k_0=k_{N+1}=0 $. We define:
$$x_{\frac{1}{2}}=0,\,\mbox{ for } i\in \{1,...,M\}:\,\,x_\II=x_\Ii+h_i,$$  
$$y_{\frac{1}{2}}=0,\,\mbox{ for } j\in \{1,...,N\}:\,\,y_\JJ=y_\Jj+k_j,$$  
(So that $x_{M+1/2}=y_{N+1/2}=1$), and 
$$Kx_i=[x_\Ii,x_\II],\,Ky_j=[y_\Jj,y_\JJ],\,K_{ij}=Kx_i\times Ky_j.$$
Let $ (x_i)_{i=0}^{M+1}$ and $ (y_j)_{j=0}^{N+1}$ be points such that
$$
x_{\Ii } < x_i<x_{\Ii }, \mbox{ for } i=1,...,M;\,x_0=0,\,x_{M+1}=1,
$$
$$
y_{\Jj } < y_j<y_{\Jj }, \mbox{ for } j=1,...,N;\,y_0=0,\,y_{N+1}=1,
$$
and let $ x_{i,j}=(x_i,y_j) \mbox{ for } i=1,...,M \mbox{ and }j=1,...,N$. Set 
$$h_i^{-}=x_i-x_\Ii,\,h_i^{+}=x_\II-x_i,\mbox{ for } i=1,...,M,\,\HI =x_{i+1}-x_i,\mbox{ for } i=0,...,M,$$ 
$$k_j^{-}=y_j-y_\Jj,\,k_j^{+}=y_\JJ-y_j,\mbox{ for } j=1,...,N,\,\KJ =y_{j+1}-y_j,\mbox{ for } j=0,...,N.$$
Assuming that $h_0^{+}=h_{M+1}^{-}=k_0^{+}=k_{N+1}^{-}=0$ and considering the
mesh size $ h=\max\{(h_i, i=1,...,M),\,( k_j=1,...,N)\}$.}} 
\begin{definition}
Let $\Xx(\tau)$ be the set of functions from $\Omega$ to $\R$ piecewise 
constant over each $K_{ij}$.

For $w \in \Xx(\tau)$, we define the discretes $H^1_0$-norm and $L^2$-norm
respectively
\begin{equation}
\label{NormH1}
\|w\|_{1,\tau}=\lf \SM k_j \frac{(w_{i+1,j}-w_{ij})^2}{\HI }+\SM h_i\frac{(w_{i,j+1}-w_{ij})^2}{\KJ }\rf ^{\frac{1}{2}} 
\end{equation}
\begin{equation}
\label{NormL2}
\|w\|_{L^2}=\lf \SM h_ik_j w_{ij}^2\rf ^{1/2}
\end{equation}
Let $w=(w_{ij})_{0\leq i\leq M+1,0\leq j\leq N+1}$ and $\Dh $ be the the 
following discrete opertor 
\begin{equation}
\label{DH}
(\Dh w)_{ij}=-k_j\lf \frac{w_{i+1,j}-w_{ij}}{\HI }-\frac{w_{i,j}-w_{i-1,j}}{\Hi }\rf-
h_i\lf \frac{w_{i,j+1}-w_{ij}}{\KJ }-\frac{w_{i,j}-w_{i,j-1}}{k_{\Jj }}\rf.
\end{equation}
For a continuous function $g\in C(\Omega)$, we introduce a similar definition:
\begin{equation}
\label{DHH}
(\Dhh g)_{ij}=-k_j(\frac{g(x_{i+1,j})-g(x_{ij})}{\HI }-\frac{g(x_{i,j}-g(x_{i-1,j})}{\Hi })-
h_i( \frac{g(x_{i,j+1})-g(x_{ij})}{\KJ }-\frac{g(x_{i,j})-g(x_{i,j-1})}{k_{\Jj }}).
\end{equation}
\end{definition}
To simlify the notations, $\sum_{ij}^1$ denotes
$\sum_{\stackrel{i=\overline{1,M}}{j=\overline{1,N}}}$ and $\sum_{ij}$ denotes
$\sum_{\stackrel{i=\overline{0,M}}{j=\overline{0,N}}}$\\
Integrating the equation (\ref{pb2d}) over each finite volume $K_{ij}$, to get
\begin{equation}
\label{Int}
-\int_{Ky_j}\lf u_x(x_\II,y)-u_x(x_\Ii,y)\rf dy-\int_{Kx_i}\lf u_y(x,y_\JJ)
-u_y(x,y_\Jj)\rf dx=\int_{K_{ij}}fdxdy.
\end{equation}
Taking the first term in left hand side (l.h.s) of (\ref{Int})
$$
\int_{Ky_j}\lf u_x(x_{\II },y)-u_x(x_{\Ii},y)\rf dy=-k_j
\lf u_x(x_{\II },y_j)-u_x(x_{\Ii },y_j)\rf -S_{\IIj}+S_{\Iij},
$$
where 
\begin{equation}
\label{SI}
S_{\IIj}=-\int_{y_{\Jj }}^{y_{\JJ }}(y-y_j)u_{xy}(x_{\II },\widehat y_j)dy,
\end{equation}
and $\widehat y_j$ is some point lies between $y_j$ and $y$.\\
Using, again Taylor's formula, yields  
\begin{equation}
\label{UxI}
u_x(x_{\II },y_j)=\frac{u(x_{i+1},y_j)-u(x_i,y_j)}{\HI }+R_{\IIj}.
\end{equation}
Thus the following estimates hold
\begin{equation}
\label{Est_SI}
|S_{\IIj}| \leq ck_j^2 |u_{xy}|_{\infty,\bar \Omega},
\mbox{ and }
|R_{\IIj}| \leq c\HI |u_{2x}|_{\infty,\bar \Omega}.
\end{equation}
By the same way, we can get  
\begin{eqnarray}
-\int_{Kx_i}\lf u_y(x,y_{\JJ })-u_y(x,y_{\Jj })\rf dx
&=&-h_i\lf \frac{u(x_i,y_{j+1})-u(x_i,y_j)}{\KJ }-
\frac{u(x_i,y_j)-u(x_i,y_{j-1})}{\KJ }\rf \nonumber \\
&&-h_i\lf R_{\iJJ}-R_{\iJj}\rf -S_{\iJJ}+S_{\iJj },
\end{eqnarray}
where 
\begin{equation}
\label{Est_SiJJ}
|S_{\iJJ}| \leq ch_i^2 |u_{xy}|_{\infty,\bar \Omega},
\mbox{ and  }
|R_{\iJJ}| \leq c\KJ |u_{2y}|_{\infty,\bar \Omega}.
\end{equation}
Therefore, the equation (\ref{pb2d}) becomes after integration as follows
$$
(\Dhh u)_{ij}= \int_{K_{ij}}f\,dxdy +k_j(R_{\IIj}-R_{\Iij }) + h_i( R_{\iJJ}-R_{\iJj })+ S_{\IIj}-S_{\Iij }+S_{\iJJ}-S_{\iJj }.
$$
The basic finite volume solution $u^h=(u_{ij})_{i=0,...,M+1, j=0,...,N+1}$ is 
defined by 
\begin{equation}
\label{Don_bor}
u_{0j}=u_{M+1,j}=u_{i0}=u_{i,N+1}=0,
\end{equation}
and for $(i,j) \in \{1,...,M\} \times \{1,...,N\}$, we have
\begin{equation}
\label{Sch_vf}
(\Dh u^h)_{ij}= \int_{K_{ij}}f\,dxdy
\end{equation}
The existence, uniqueness of the solution $u^h$, the analysis of the order of 
the convergence can  be justified as done for 1D case (see \cite{Eymard}). More
precisely, we have the following theorem
\begin{theorem}
\label{TH_BOOK}
(\cite{Eymard})
 If the solution $u$ of the equation (\ref{pb2d}) belonging to 
$C^2(\bar \Omega)$ and $ f \in C^2(\bar \Omega)$. Then the approximate solution
$u^h=(u_{ij})$ defined by the boundary condition (\ref{Don_bor}) and the
discrete equation (\ref{Sch_vf}), satisfies the following estimates
\begin{equation}
\label{H0VF}
\|e\|_{1,\tau} \leq ch \| u \|_{2,\infty,\bar \Omega},
\end{equation}
\begin{equation}
\label{L2VF}
\|e\|_{L^2} \leq ch \| u \|_{2,\infty,\bar \Omega},
\end{equation}
\begin{equation}
\label{L2MOD}
\lf \SMN k_j\HI e_{ij}^2 \rf ^ \frac{1}{2} \leq ch \| u \Vert_{2,\infty,\bar \Omega}.
\end{equation}
where $ e_{ij}=u(x_i,y_j)-u_{ij} $ for $(i,j)\in \{1,...,M\}\times \{1,...,N\}$
and vanishes elsewhere.
\end{theorem} 
\subsection{First Correction} 
In this section, we assume more regularity for the solution $u$, i.e. 
$u \in C^4(\bar \Omega)$.\\
Looking, again, at the equation (\ref{Int}), to simplify the notation, let  
$$
g_i(y)=u_x(x_{\II },y)-u_x(x_{\Ii },y).
$$
Using Taylor's formula, we get:
\begin{equation}
\label{Int_g}
 -\int_{y_{\Jj }}^{y_{\JJ }}g_i(y)dy=- k_jg_i(y_j)-\frac{{k_j^{+}}^2-{k_j^{-}}^2}{2}(g_i)_y(y_j)-T_{ij},
\end{equation}
where 
\begin{equation}
\label{Est_TIJp}
\begin{array}{ccl}
\vert T_{ij} \vert &\leq& ck_j^3\vert (g_i)_{2y} \vert_{\infty,K_j}\\
&\leq& ck_j^3h_i\vert u_{2x,2y} \vert_{\infty,\bar \Omega}.
\end{array}
\end{equation}
Using again, Taylor's formula, yields
\begin{equation}
\label{giy}
\begin{array}{ccl}
({k_j^{+}}^2-{k_j^{-}}^2)(g_i)_y(y_j)&=&
({k_j^{+}}^2-{k_j^{-}}^2)(u_{xy}(x_{\II },y_j)-u_{xy}(x_{\Ii },y_j))\\
&&\\
&=&{k_j^{+}}^2h_iu_{2x,y}(x_i,y_j)-{k_j^{-}}^2h_iu_{2x,y}(x_i,y_{j-1})+U_{ij}^1\\
&&+({k_j^{+}}^2-{k_j^{-}}^2)S_{ij}^1-({k_j^{+}}^2-{k_j^{-}}^2)S_{ij}^2
\end{array}
\end{equation}
where
\begin{equation}
\label{Est_SIJ1}
\vert S_{ij}^1 \vert \leq c{h_i^{+}}^2\vert u_{3x,y} \vert_{\infty,\bar \Omega}
\end{equation}
\begin{equation}
\label{Est_SIJp}
\vert S_{ij}^2 \vert \leq c{h_i^{-}}^2\vert u_{3x,y} \vert_{\infty,\bar \Omega}
\end{equation}
\begin{equation}
\label{Est_UIJ1}
\vert U_{ij}^1 \vert \leq c{k_j^{-}}^2h_ik_{\Jj }\vert u_{2x,2y} \vert_{\infty,\bar \Omega}
\end{equation}
By substituting $g_i$ by its value in (\ref{Int_g}) and by using equality (\ref{giy}), we get
\begin{equation}
\label{Est_Int_g}
\begin{array}{rl}
\dsp{-\int_{y_{\Jj }}^{y_{\JJ }}\lf u_x(x_{\II },y)-u_x(x_{\Ii },y)\rf dy}&=-k_j\lf u_x(x_{\II },y_j)-u_x(x_{\Ii },y_j)\rf \\
&\dsp{-h_i \frac{{k_j^{+}}^2}{2}u_{2x,y}(x_i,y_j)+h_i \frac{{k_j^{-}}^2}{2}u_{2x,y}(x_i,y_{j-1})} \\
&\\
&\dsp{-\frac{{k_j^{+}}^2-{k_j^{-}}^2}{2}S_{ij}^1-U_{ij}^1+\frac{{k_j^{+}}^2-{k_j^{-}}^2}{2}S_{ij}^2 -T_{ij}},
\end{array}
\end{equation}
but, in the other hand
\begin{equation}
\label{UxI+1}
\frac{u(x_{i+1},y_j)-u(x_i,y_j)}{\HI }=u_x(x_{\II },y_j)+
\frac{h_{i+1}^{-}-h_i^{+}}{2}u_{2x}(x_{\II },y_j)+R_{\IIj}^2,
\end{equation}
where
\begin{equation}
\label{Est_RIJ2}
\vert R_{\IIj}^2 \vert \leq c \HI ^2 
\vert u_{3x} \vert_{\infty,\bar \Omega}
\end{equation}
Combining (\ref{Est_Int_g}) and (\ref{UxI+1}) yields that
\begin{eqnarray}
\label{2deq101}
-\int_{y_{\Jj }}^{y_{\JJ }}\lf u_x(x_{\II },y)-u_x(x_{\Ii },y)\rf dy&=&
-k_j\lf \frac{u(x_{i+1},y_j)-u(x_i,y_j)}{\HI}-
\frac{u(x_i,y_j)-u(x_{i-1},y_j)}{\Hi }\rf \nonumber\\ 
&+&k_j\lf \frac{h_{i+1}^{-}-h_i^{+}}{2}u_{2x}(x_{i+ \frac{1}{2}},y_j)-
\frac{h_i^{-}-h_{i-1}^{+}}{2}u_{2x}(x_{i- \frac{1}{2}},y_j)\rf  \nonumber \\ 
&-&h_i \frac{{k_j^{+}}^2}{2}u_{2x,y}(x_i,y_j)+
h_i \frac{{k_j^{-}}^2}{2}u_{2x,y}(x_i,y_{j-1})
+k_jR_{\IIj}^2\nonumber \\         
&-& k_jR_{\Iij }^2-\frac{{k_j^{+}}^2-{k_j^{-}}^2}{2}S_{ij}^1+
\frac{{k_j^{+}}^2-{k_j^{-}}^2}{2}S_{ij}^2-U_{ij}^1+T_{ij} 
\end{eqnarray}
and by the same way, we can find similar expansion for the second term in the
l.h.s of (\ref{Int}).
\begin{eqnarray}
\label{2deq102}
-\int_{x_{\Ii }}^{x_{\II }}\lf u_y(x,y_{\JJ })-
u_y(x,y_{\Jj })\rf dy&=&-h_i\lf\frac{u(x_i,y_{j+1})-u(x_i,y_j)}{\KJ }-
\frac{u(x_i,y_j)-u(x_i,y_{j-1})}{k_{\Jj }}\rf \nonumber          \\ 
&+&h_i\lf\frac{k_{j+1}^{-}-k_j^{+}}{2}u_{2y}(x_i,y_{j+ \frac{1}{2}})-
\frac{k_j^{-}-k_{j-1}^{+}}{2}u_{2y}(x_i,y_{j- \frac{1}{2}})\rf   \nonumber \\ 
&-&k_j \frac{{h_i^{+}}^2}{2}u_{x,2y}(x_i,y_j)+k_j \frac{{h_i^{-}}^2}{2}u_{x,2y}(x_{i-1},y_j)
+h_iR_{\iJJ}^2\nonumber \\
&-& h_iR_{\iJj }^2-\frac{{h_i^{+}}^2-{h_i^{-}}^2}{2}L_{ij}^1+
\frac{{h_i^{+}}^2-{h_i^{-}}^2}{2}L_{ij}^2-U_{ij}^2-H_{ij}
\end{eqnarray}
where 
\begin{equation}
\label{Est_RJ2}
\vert R_{\iJJ}^2 \vert \leq c \KJ ^2 
\vert u_{3y} \vert_{\infty,\bar \Omega}
\end{equation}
\begin{equation}
\label{Est_LIJ1}
\vert L_{ij}^1 \vert \leq c{k_j^{+}}^2\vert u_{x,3y} \vert_{\infty,\bar
\Omega}
\end{equation}
\begin{equation}
\label{Est_LIJ2}
\vert L_{ij}^2 \vert \leq c{k_j^{-}}^2\vert u_{x,3y} \vert_{\infty,\bar
\Omega}
\end{equation}
\begin{equation}
\label{Est_HIJ}
\vert H_{ij} \vert \leq ch_i^3 k_j\vert u_{2x,2y} \vert_{\infty,\bar \Omega}\\
\end{equation}
\begin{equation}
\label{Est_UIJ2}
\vert U_{ij}^2 \vert \leq c{h_i^{-}}^2k_j\Hi \vert u_{2x,2y} \vert_{\infty,\bar \Omega}
\end{equation}
Equalities (\ref{2deq101}) and (\ref{2deq102}) combined with (\ref{Int}) yields
\begin{eqnarray}
\label{EQ112}
(\Dhh u)_{ij}&=&\int_{K_{ij}}fdxdy
-k_j\frac{h_{i+1}^{-}-h_i^{+}}{2}u_{2x}(x_{i+ \frac{1}{2}},y_j)+
k_j\frac{h_i^{-}-h_{i-1}^{+}}{2}u_{2x}(x_{i- \frac{1}{2}},y_j)\nonumber \\ 
&&-h_i\frac{k_{j+1}^{-}-k_j^{+}}{2}u_{2y}(x_i,y_{j+ \frac{1}{2}})+
h_i\frac{k_j^{-}-k_{j-1}^{+}}{2}u_{2y}(x_i,y_{j- \frac{1}{2}})\nonumber \\ 
&&+h_i \frac{{k_j^{+}}^2}{2}u_{2x,y}(x_i,y_j)-
h_i \frac{{k_j^{-}}^2}{2}u_{2x,y}(x_i,y_{j-1})
+k_j \frac{{h_i^{+}}^2}{2}u_{x,2y}(x_i,y_j)-
k_j \frac{{h_i^{-}}^2}{2}u_{x,2y}(x_{i-1},y_j)\nonumber \\
&&-k_jR_{\IIj}^2+ k_jR_{\Iij }^2         
+\frac{{k_j^{+}}^2-{k_j^{-}}^2}{2}S_{ij}^1-
\frac{{k_j^{+}}^2-{k_j^{-}}^2}{2}S_{ij}^2+T_{ij}\nonumber \\
&&-h_iR_{\iJJ}^2+ h_iR_{\iJj }^2         
+\frac{{h_i^{+}}^2-{h_i^{-}}^2}{2}L_{ij}^1-
\frac{{h_i^{+}}^2-{h_i^{-}}^2}{2}L_{ij}^2+H_{ij}
+U_{ij}^1+U_{ij}^2.
\end{eqnarray}
Let
$\Gamma_1=[0,1]\times\{0\},\Gamma_2=\{1\}\times[0,1],\Gamma_3=[0,1]\times\{1\},
\Gamma_4=\{0\}\times[0,1]$ be the partial boundaries. The second derivative 
$u_{2x}$ of $u$ is the solution of the problem
\begin{equation}\label{eqn} 
\left\lbrace
\begin{array}{l}
-\Delta v=f_{2x}, \mbox{on } \Omega =(0,1)^2\\
v|_ {\Gamma_1}=v|_ {\Gamma_3}=0,\\
v|_ {\Gamma_2}=-f|_ {\Gamma_2},\\
v|_ {\Gamma_4}=-f|_ {\Gamma_4}
\end {array}
\right.
\end{equation}
Let $v^h=(v_{ij})$ the finite volume approximation of $u_{2x}$, then
\begin{equation}
\label{Don_Bor_v} 
\left\lbrace
\begin{array}{l}
v_{i0}=v_{i,N+1}=0,\\
v_{0j}=-f(0,y_j),\\
v_{M+1,j}=-f(1,y_j)
\end {array}
\right.
\end{equation}
and for $(i,j) \in \{1,...,M\} \times \{1,...,N\}$, we have
\begin{equation}
\label{Sch_v}
(\Dh v^h)_{ij}=\int_{K_{ij}}f_{2x}dxdy.
\end{equation}
\begin{remarque}
We remark that, we have used the same matrix, that used to compute $u^h$, to 
compute an optimal approximation for $u_{2x}$. This implies, in turn, that we 
have the same order of the convergence as for the solution $u$ provided that 
$u$ belonging to $C^4(\bar \Omega)$, i.e. theorem \ref{TH_BOOK} holds also for $u_{2x}$ 
instead of $u$ and $v^h$ instead of $u^h$.
\end{remarque}
 Therefore
\begin{lemma}
\label{lemma8}
Let $u$ be the solution of the equation (\ref{pb2d}). If 
$u \in C^4(\bar \Omega)$. Then the approximation $v^h$ defined by the boundary 
conditions (\ref{Don_Bor_v}) and the equation (\ref{Sch_v}) satifies the 
following estimates
\begin{equation}
\label{L8E1}
\lf \SMN h_ik_j (u_{2x}(x_{\II },y_j)-v_{ij})^2 \rf^ \frac{1}{2} \leq ch \Vert u \Vert_{4,\infty,\bar \Omega}
\end{equation}
\begin{equation}
\label{L8E2}
\lf \SMN h_i\KJ (u_{2x,y}(x_i,y_j)-\frac{v_{i,j+1}-v_{ij}}{\KJ })^2\rf^ \frac{1}{2} \leq ch \Vert u \Vert_{4,\infty,\bar \Omega}
\end{equation}
\begin{equation}
\label{L8E3}
\lf \SMN \HI k_j (u_{3x}(x_i,y_j)-\frac{v_{i+1,j}-v_{ij}}{\HI })^2\rf^ \frac{1}{2} \leq ch \Vert u \Vert_{4,\infty,\bar \Omega}
\end{equation}
\end{lemma} 
{\bf{Proof}}\\
1-Using triangular inequality combined with estimate (\ref{L2VF}) of theorem 
\ref{TH_BOOK} yield 
\begin{eqnarray*}
\lf\SMN h_ik_j
(u_{2x}(x_{\II },y_j)-v_{ij})^2\rf^ \frac{1}{2} &\leq& 
\lf\SMN h_ik_j
(u_{2x}(x_{\II },y_j)-u_{2x}(x_i,y_j)^2\rf^ \frac{1}{2} \\
&+&\lf \SMN h_ik_j
(u_{2x}(x_i,y_j)-v_{ij})^2\rf^ \frac{1}{2} \\
&\leq& ch \,(\Vert u_{3x} \Vert_\infty+\Vert u \Vert_{4,\infty,\bar \Omega})\\
&\leq& ch \Vert u \Vert_{4,\infty,\bar \Omega}.
\end{eqnarray*}
2-Estimates (\ref{L8E2}) and (\ref{L8E3}) can be proven by the same way, i.e. 
we use triangular inequality and theorem \ref{TH_BOOK}.$\Box$\\
Let $w^h=(w_{ij})$ be a discrete function, to simplify the natation, we define the following discrete operators 
$$
\partial_1^x w_{ij}=\frac{w_{i+1,j}-w_{ij}}{\HI }, 
\mbox{  and   }\partial_1^y w_{ij}=\frac{w_{i,j+1}-w_{ij}}{\KJ }.
$$ 

Now, we are able to define a new approximation $u^h_1=(u_{ij}^1)$, called
correction, of order $O(h^2)$ (as we will see) defined by the boundary 
conditions (\ref{Don_bor}) and the following discrete equation for 
$ (i,j) \in \{1,...,M\}\times \{1,...,N\}$   
\begin{equation}
\label{Sch_Cor}
(\Dh u^h_1)_{ij}=\int_{K_{ij}}f\,dx\,dy +\lf\gamma_{\IIj}-\gamma_{\Iij }\rf+
\lf \gamma_{\iJJ}-\gamma_{\iJj }\rf+\delta_{ij},
\end{equation}
where
\begin{equation}
\gamma_{\IIj}=-k_j\frac{h_{i+1}^{-}-h_i^{+}}{2}v_{ij},
\end{equation}
\begin{equation}
\gamma_{\iJJ}=-h_i\frac{k_{j+1}^{-}-k_j^{+}}{2}(-f(x_i,y_{\JJ })-v_{ij}),
\end{equation} 
\begin{equation}
\left\{
\begin{array}{ccl}
\delta_{ij}&=&\dsp{h_i \frac{{k_j^{+}}^2}{2}  \partial_1^y v_{ij} -h_i \frac{{k_j^{-}}^2}{2} \partial_1^y v_{i,j-1}
+k_j \frac{{h_i^{+}}^2}{2}(-f_x(x_i,y_j)-\partial_1^x v_{ij})} \\  
&&\\
&&\dsp{-k_j \frac {{h_i^{-}}^2}{2}(-f_x(x_{i-1},y_j)-\partial_1^x v_{i-1,j}).}
\end{array} 
\right.
\end{equation}
\begin{remarque} 
As you can see that, we have used only tha approximation of $u_{2x}$ and the 
fact that $u_{2y}=-f-u_{2x}$ to approximate all the higher pointwise 
derivatives in the r.h.s of (\ref{EQ112}).
\end{remarque}

To analyse the error, let $e_{ij}^1=u_{ij}^1-u(x_i,y_j)$ be the error in the
first correction, by substracting (\ref{EQ112}) from (\ref{Sch_Cor}), we get
\begin{eqnarray}
\label{EIJ1}
(\Dh e^1)_{ij}=\bar \gamma_{\IIj}-\bar \gamma_{\Iij }
+\bar \gamma_{\iJJ}-\bar \gamma_{\iJj }+\bar \delta_{ij},
\end{eqnarray} 
where
\begin{equation}
\label{GIIJ}
\bar \gamma_{\IIj }=-k_j\frac{h_{i+1}^{-}-h_i^{+}}{2}
(v_{ij}-u_{2x}(x_{\II },y_j))+k_jR_{\IIj}^2,
\end{equation}
\begin{equation}
\label{GIJJ}
\bar\gamma_{\iJJ}=-h_i\frac{k_{j+1}^{-}-k_j^{+}}{2}
(u_{2x}(x_i,y_{\JJ })-v_{ij})+h_iR_{\iJJ}^2,
\end{equation} 
\begin{eqnarray}
\bar \delta_{ij}=h_i \frac{{k_j^{+}}^2}{2} ( \partial_1^y v_{ij}-u_{2x,y}(x_{i,j}))-h_i \frac{{k_j^{-}}^2}{2}( \partial_1^y v_{i,j-1}-u_{2x,y}(x_{i,j-1}))-h_i \frac{{k_j^{-}}^2}{2}( \partial_1^y v_{i,j-1}-u_{2x,y}(x_{i,j-1}))\nonumber \\
+k_j \frac{{h_i^{+}}^2}{2}(u_{3x}(x_i,y_j)-\partial_1^x v_{ij})-k_j \frac {{h_i^{-}}^2}{2}(u_{3x}(x_{i-1},y_j)-\partial_1^x v_{i-1,j})\nonumber \\ 
-\frac{{k_j^{+}}^2-{k_j^{-}}^2}{2}S_{ij}^1+\frac{{k_j^{+}}^2-{k_j^{-}}^2}{2}S_{ij}^2\frac{{h_i^{+}}^2-{h_i^{-}}^2}{2}L_{ij}^1+\frac{{h_i^{+}}^2-{h_i^{-}}^2}{2}L_{ij}^2\nonumber \\ 
-T_{ij}-H_{ij}
-U_{ij}^1-U_{ij}^2
\end{eqnarray} 
\subsection{Convergence Order of the First Correction}

Multiplying both sides of (\ref{EIJ1}) by $e_{ij}^1$ and summing over $(i,j) \in \{1,...,M\}
\times \{1,...,N\}$, to get
\begin{eqnarray}
\label{SEIJ1}
-\sMN k_j\frac{e_{i+1,j}^1-e_{ij}^1}{\HI }e_{ij}^1+\sMN k_j\frac{e_{ij}^1-e_{i-1,j}^1}{\Hi }e_{ij}^1
-\sMN h_i\frac{e_{i,j+1}^1-e_{ij}^1}{\KJ }e_{ij}^1+\sMN h_i\frac{e_{ij}^1-e_{i,j-1}^1}{k_{\Jj }}e_{ij}\nonumber \\
=\sMN \bar \gamma_{\IIj}e_{ij}^1-\sMN \bar \gamma_{\Iij }e_{ij}^1+\sMN \bar \gamma_{\iJJ}e_{ij}^1-
\sMN \bar \gamma_{\iJj }e_{ij}^1+\sMN \bar \delta_{ij}e_{ij}^1
\end{eqnarray}
Reordering equation (\ref{SEIJ1}) and using the fact that $e_{ij}^1$ vanishs on the boundary mesh
points, to get
\begin{equation}
\label{sEIJ1}
\|e^1\|_{1,\tau}^2= -\SMN \bar \gamma_{\IIj}(e_{i+1,j}^1-e_{ij}^1) -\SMN \bar \gamma_{\iJJ}(e_{i,j+1}^1-e_{ij}^1)+\sMN \bar \delta_{ij}e_{ij}^1
\end{equation}
We should now estimate each term in the r.h.s of (\ref{sEIJ1}). The first and 
the second term can be handled by the same way. Hence, it is suffices to 
estimate the first and the last ones. Using Cauchy-Schwarz inequality yields 
that
\begin{eqnarray}
\vert \SMN \bar \gamma_{\IIj}(e_{i+1,j}^1-e_{ij}^1) \vert &\leq& (\SMN \frac{\bar \gamma_{\IIj}^2\HI }{k_j})^ \frac{1}{2}(
\SMN k_j\frac{(e_{i+1,j}^1-e_{ij}^1)^2}{\HI })^ \frac{1}{2} \nonumber \\
&\leq&  (\SMN \frac{\bar \gamma_{\IIj}^2\HI }{k_j})^ \frac{1}{2} 
\|e^1\|_{1,\tau} 
\end{eqnarray}
Using triangular inequality, in order to get
\begin{eqnarray*}
(\SMN \frac{\bar  \gamma_{\IIj}^2\HI }{k_j})^ \frac{1}{2} \leq 
\lf \SMN \frac{k_j^2(h_{i+1}^{-}-h_i^{+})^2}{4k_j}\HI (v_{ij}-u_{2x}(x_{\II },y_j))^2\rf^ \frac{1}{2}
+\lf\SMN \frac{k_j^2(R_{\II }^2)^2}{k_j}\HI \rf^ \frac{1}{2} 
\end{eqnarray*}
Using estimate (\ref{L8E1}) of lemma \ref{lemma8} (combined with estimate (\ref{L2MOD}) of 
theorem \ref{TH_BOOK} ) and inequality (\ref{Est_RIJ2}), to get
\begin{eqnarray}
\label{Est_hk}
\lf\SMN \frac{\bar \gamma_{\IIj}^2\HI }{k_j}\rf^ \frac{1}{2} &\leq& 
c\lf h^2 \Vert u \Vert_{4,\infty,\bar \Omega}+h^2\Vert u \Vert_{3,\infty,\bar \Omega}\rf \nonumber \\
&\leq& ch^2 \Vert u \Vert_{4,\infty,\bar \Omega}
\end{eqnarray}
Comming back now to the last term in the r.h.s of (\ref{sEIJ1}), using triangular inequality
yields that
\begin{eqnarray}
\label{SBAR}
\vert\sMN \bar \delta_{ij}e_{ij}^1\vert &\leq&
\sMN  h_i \frac{{k_j^{+}}^2}{2} \vert \partial_1^y v_{ij}-u_{2x,y}(x_i,y_j)\vert \vert e_{ij}^1\vert + 
\sMN  h_i \frac{{k_j^{-}}^2}{2} \vert \partial_1^y v_{i,j-1}-u_{2x,y}(x_i,y_{j-1})\vert \vert e_{ij}^1\vert \nonumber \\
&+& \sMN k_j \frac{{h_i^{+}}^2}{2} \vert u_{3x}(x_i,y_j)-\partial_1^x v_{ij}\vert \vert e_{ij}^1\vert +
    \sMN k_j \frac {{h_i^{-}}^2}{2} \vert u_{3x}(x_{i-1},y_j)-\partial_1^x v_{i-1,j}\vert \vert e_{ij}^1\vert\nonumber \\ 
&+& \sMN \frac{\vert{k_j^{+}}^2-{k_j^{-}}^2\vert}{2} \vert S_{ij}^1 \vert  \vert
e_{ij}^1\vert +\sMN \frac{\vert{k_j^{+}}^2-{k_j^{-}}^2\vert}{2}
\vert S_{ij}^2 \vert \vert e_{ij}^1\vert +\sMN  \vert  T_{ij} \vert \vert e_{ij}^1\vert\nonumber \\ 
&+& \sMN \frac{\vert{h_i^{+}}^2-{h_i^{-}}^2\vert}{2} \vert L_{ij}^1 \vert \vert e_{ij}^1\vert+
\sMN \frac{\vert{h_i^{+}}^2-{h_i^{-}}^2\vert}{2} \vert L_{ij}^2 \vert \vert e_{ij}^1\vert
+\sMN  \vert H_{ij}\vert \vert e_{ij}^1\vert \nonumber \\ 
&+&\sMN  \vert U_{ij}^1\vert \vert e_{ij}^1\vert+\sMN  \vert U_{ij}^2 \vert  \vert e_{ij}^1\vert
\end{eqnarray}
Begininig by the first term in the r.h.s of (\ref{SBAR}), using the Cauchy-Schwars inequality and the
estimate (\ref{L8E2}) of lemma \ref{lemma8}
\begin{eqnarray}
\sMN h_i \frac{{k_j^{+}}^2}{2} \vert \partial_1^y v_{ij}-u_{2x,y}(x_i,y_j)\vert \vert e_{ij}^1\vert&\leq& \sMN h_i \frac{{\KJ }^2}{2} \vert \partial_1^y v_{ij}-u_{2x,y}(x_i,y_j)\vert \vert e_{ij}^1\vert \nonumber \\ 
&\leq&  ch \lf \sMN h_i \KJ ( \partial_1^y v_{ij}-u_{2x,y}(x_i,y_j))^2\rf^ \frac{1}{2} 
 \lf \sMN h_i \KJ ( e_{ij}^1)^2\rf^ \frac{1}{2} \nonumber \\
&\leq&  ch^2 \NNSU \lf \sMN h_i \KJ  (e_{ij}^1)^2\rf^ \frac{1}{2}    
\end{eqnarray}
This implies that
\begin{equation}
\sMN h_i \frac{{k_j^{+}}^2}{2} \vert \partial_1^y v_{ij}-u_{2x,y}(x_i,y_j)\vert \vert e_{ij}^1\vert \leq ch^2 \|e^1\|_{1,\tau}
\Vert u \Vert_{4,\infty, \bar \Omega}
\end{equation}
and by the same way, we can handle the second term in the r.h.s of (\ref{SBAR}). Indeed
\begin{eqnarray*}
\sMN h_i \frac{{k_j^{-}}^2}{2}\vert \partial_1^y v_{i,j-1}-u_{2x,y}(x_i,y_{j-1})\vert \vert e_{ij}^1\vert &\leq&
ch\sMN h_ik_j^{-} \vert \partial_1^y v_{i,j-1}-u_{2x,y}(x_i,y_{j-1})\vert \vert e_{ij}^1\vert \nonumber \\ 
&\leq&  ch (\sMN h_i k_j^{-}( e_{ij}^1)^2)^ \frac{1}{2}  (\sMN h_i
k_j^{-}( \partial_1^y v_{i,j-1}-u_{2x,y}(x_i,y_{j-1}))^2)^ \frac{1}{2} 
\nonumber \\
&\leq& ch \|e^1\|_{L^2}\lf \sMN h_i k_{\Jj }( \partial_1^y v_{i,j-1}-u_{2x,y}(x_i,y_{j-1}))^2\rf^ \frac{1}{2} \nonumber \\  
&\leq&  ch \|e^1\|_{L^2} \lf \sum_{\stackrel{i=\overline{1,M}}{j=\overline{0,N}}}h_i\KJ ( \partial_1^y v_{ij}-u_{2x,y}(x_i,y_j))^2\rf^ \frac{1}{2}. 
\end{eqnarray*}   
This with estimate (\ref{L8E2}) imply that
\begin{equation}
\sMN h_i \frac{{k_j^{-}}^2}{2}
\vert \partial_1^y v_{i,j-1}-u_{2x,y}(x_i,y_{j-1})\vert \vert e_{ij}^1\vert 
\leq ch^2 \|e^1\|_{1,\tau}\Vert u \Vert_{4,\infty, \bar \Omega}
\end{equation} 
Taking, now, a look at the other kind of terms. Using estimate (\ref{Est_SIJ1}) to get
\begin{eqnarray*}
\sMN \frac{\vert{k_j^{+}}^2-{k_j^{-}}^2\vert}{2}
\vert S_{ij}^1 \vert \vert e_{ij}^1\vert  &\leq&
 c \Vert u \Vert_{4,\infty, \bar \Omega}\sMN k_j^2h_i^2 \vert e_{ij}^1\vert
  \\
 &\leq& ch^2 \Vert u \Vert_{4,\infty, \bar \Omega}
 \sMN k_jh_i \vert e_{ij}^1\vert\\
&\leq& ch^2 \Vert u \Vert_{4,\infty, \bar \Omega}
\lf\sMN k_jh_i\rf^\frac{1}{2}
  \lf\sMN k_jh_i(e_{ij}^1)^2\rf^\frac{1}{2}.
\end{eqnarray*}
Hence
\begin{equation}
\label{SIJ1}
\sMN \frac{\vert{k_j^{+}}^2-{k_j^{-}}^2\vert}{2}
\vert S_{ij}^1 \vert \vert e_{ij}^1\vert  \leq ch^2 \NNSU \|e^1\|_{1,\tau}.
\end{equation} 
By the same way, we can find the same estimate for terms corresponding to 
$ S_{ij}^2 $, $ L_{ij}^1 $ , $L_{ij}^2 $, $T_{ij}$, $U_{ij}^1$ and $U_{ij}^2$ .\\
Inequalities (\ref{SBAR})-(\ref{SIJ1}) yield that
\begin{equation}
\label{SDIJ}
\vert\sMN \bar \delta_{ij}e_{ij}^1\vert
\leq  ch^2 \NNSU \|e^1\|_{1,\tau}.
\end{equation}
Combining equality (\ref{sEIJ1}) with inequalities (\ref{Est_hk}) and 
(\ref{SDIJ}) yields
\begin{theorem}
\label{TH6}
If the solution $u$ of the equation (\ref{pb2d}) belonging to 
$C^4(\bar \Omega)$ . Let $u^h$ be the basic finite volume solution of boundary 
condition (\ref{Don_bor}) and discrete equation (\ref{Sch_vf}). Then the finite
volume approximation $u^h_1=(u_{ij}^1)$ defined by the boundary condition 
(\ref{Don_bor}) and the discrete equation (\ref{Sch_Cor}), satisfies the 
following $O(h)$ improvement in $H^1_0$-norm
\begin{equation}
\label{T6E1}
 \|e^1\|_{1,\tau}\leq ch^2 \NNSU,
\end{equation}
\begin{equation}
\label{T6E2}
\lf\SMN k_j\HI (e_{ij}^1)^2\rf^ \frac{1}{2} \leq ch^2 \NNSU,
\end{equation}
\begin{equation}
\label{T6E3}
\|e^1\|_{L^2} \leq ch^2 \NNSU,
\end{equation}
where $ e_{ij}^1=u(x_i,y_j)-u_{ij}^1 $ for 
$(i,j) \in \{1,...,M\}\times \{1,...,N\}$ and vanishes elsewhere.
\end{theorem} 
\begin{remarque} Numerical results schows that the order of the convergence in
$L^2$ norm even of the basic solution $u^h$ is $O(h^2)$, i.e. the order in 
(\ref{H0VF}) and (\ref{L2VF}) is $O(h^2)$ but the coeffecient of the order in 
(\ref{T6E2}) and (\ref{T6E3}) are smaller than to those ones of (\ref{L2VF}) 
and (\ref{L2MOD}).
\end{remarque}
\subsection{The Second Correction and Higher Order of Corrections}
In this subsection, we give an idea allowing us to construct second 
correction,i.e.the order of the convergence is $O(h^3)$, this result can
be extended to construct an arbitrary correction we wish. To compute
the second correction, we use the first correction to estimate some
pointwise derivatives of the solution $u$ and the fact the equation 
(\ref{pb2d}) satisfying by $u$. It suffices to remark that, provided that at 
least $u\in C^4 (\bar \Omega)$
\begin{eqnarray*}
-\int_{K_{ij}}fdxdy&=&-k_j\lf u_x(x_{\II },y_j)-u_x(x_{\Ii },y_j)\rf
-h_i\lf u_y(x_i,y_{\JJ })-u_y(x_i,y_{\Jj })\rf\\
&&-h_i\frac{{k_j^+}^2-{k_j^-}^2}{2}u_{2x,y}(x_i,y_j)-
k_j\frac{{h_i^+}^2-{h_i^-}^2}{2}u_{x,2y}(x_i,y_j)\\
&&-\frac{{h_i^+}^2-{h_i^-}^2}{2}k_j\frac{{h_i^+}^2-{h_i^-}^2}{2}f_{xy}(x_i,y_j)
+O(k_j^4h_i)-O(h_i^4k_j)+O(k_j^2h_i^3)-O(k_j^3h_i^2).
\end{eqnarray*}
To approximate the pointwise derivative  $u_{2x,y}(x_i,y_j)$ , 
we compute the first correction to the unknown solution $u_{2x}$ (because it is
satifying the same equation that is satisfying by $u$, i.e. \ref{eqn}, this first correction is of order $O(h^2)$ in $H^1_0$, 
this means that $u_{2x,y}(x_i,y_j)$ can be approximated by an $O(h^2)$. 
By similar way, we can approximate $u_{2x,y}(x_i,y_j)$.

Comming back to $u_x(x_{\II },y_j)-u_x(x_{\Ii },y_j) $, the derivatives
will appeared here in the approximation of $u_x(x_{\II },y_j)$ and 
$u_x(x_{\Ii },y_j)$ are similars to those obtained in one dimensional space for
the second correction, and consequently, we can use the approximations just obtained to
$u_{2x}$ to approximate such derivatives.
\subsection{Some Extensions of the Results} 
So far, we have considered the Laplace model, where the second
derivative of the solution $u$ are also solutions of the same equation. In this section, we
attempt to extend results obtained, to some second order elliptic problems, where the second
derivatives of the solution $u$ are also solutions but for second member depends on the solution
$u$ itself, its derivatives and a given function. \\
The idea will be used is to approximate these terms by theirs ones corresponding in the finite
volume solution, i.e. $ u$ and derivatives of $u$ will be replaced by $u^h$ and divided difference
of $u^h$ respectively.\\
Let us consider the following model
\begin{equation}
\label{pb2dM} 
\left\lbrace
\begin{array}{l}
-\Delta u+pu=f, \mbox{on } \Omega =(0,1)^2\\
u|_ \Gamma=0.
\end {array}
\right.
\end{equation} 
where $p$ is a given function and $p \geq 0$.\\
We use the same scheme that used for Laplace model. As done above, we look, at first, for the finite volume
solution $u^h$ ( basic solution), after, we look for a convenient expansion for the error, where we try to 
approximate the derivatives of the unknown solution $u$ by using the basic solution $u^h$.
\subsection{The Finite Volume Approximation ( Basic Solution)} 
We use the same notations that used in the second section, therefore, for $f, p \in
C^1(\bar \Omega)$
\begin{eqnarray}
\label{VF_PM}
(\Dhh u)_{ij}+h_ik_jp_{ij}u(x_i,y_j)&=&\int_{K_{ij}}fdxdy
+k_j(R_{\IIj}-R_{\Iij })+h_i(R_{\iJJ}-R_{\iJj })\nonumber \\
&+&S_{\IIj}-S_{\Iij }+S_{\iJJ}-S_{\iJj }-N_{ij},
\end{eqnarray}
where $p_{ij}=p(x_i,y_j)$, $R_{\IIj}$, $R_{\iJJ}$, $S_{\iJJ}$ are defined as in
(\ref{SI}),(\ref{UxI}), (\ref{Est_SI}), (\ref{Est_SI}), (\ref{Est_SiJJ}) and 
(\ref{Est_SiJJ}). $N_{ij}$ is defined by
\begin{equation}
\label{NIJ}
N_{ij}=\int_{K_{ij}}\lf(x-x_i)\frac{\partial (pu)}{\partial x}(\widehat{a}_{ij})+
(y-y_j)\frac{\partial (pu)}{\partial y}(\widehat{a}_{ij})\rf dxdy 
\end{equation}
where  $\widehat{a}_{ij}$ is a some point in $K_{ij} $. 
Then the following estimation holds
\begin{equation}
\vert N_{ij}\vert \leq c h_ik_j(h_i+k_j) \Vert u \Vert_{1,\infty,\bar \Omega}.
\end{equation}
The basic solution $u^h=(u_{ij})_{i=0,...,M+1, j=0,...,N+1}$, which will approximate the solution $u$ of the equation
(\ref{pb2dM})  is defined by 
\begin{equation}
\label{BorM}
u_{0j}=u_{M+1,j}=u_{i0}=u_{i,N+1}=0,
\end{equation}
and for $(i,j) \in \{1,...,M\} \times \{1,...,N\}$, we have
\begin{equation}
\label{Sch_VFM}
(\Dh u^h)_{ij}+h_ik_jp_{ij}u_{ij}=\int_{K_{ij}}f\,dx\,dy.
\end{equation}
The existence and uniqueness can be done by using the same techniques in 1D 
(see \cite{Eymard}).\\
Using techniques that used in the second section yields (given in 
\cite{Eymard})
\begin{theorem}
\label{TH7}
If the solution $u$ of the equation (\ref{pb2dM}) belonging to 
$C^2(\bar \Omega)$,the coefficient $p$ belonging to $C^1(\bar \Omega) $ and 
$f \in C(\bar \Omega) $. Then the approximate solution $u^h=(u_{ij})$ defined 
by the boundary condition (\ref{BorM}) and the discrete equation 
(\ref{Sch_VFM}), satisfies the following estimates
\begin{equation}
\label{T7E1}
\|e\|_{1,\tau} \leq ch \Vert u \Vert_{2,\infty,\bar \Omega},
\end{equation}
\begin{equation}
\label{T7E2}
\lf \SMN k_j\HI e_{ij}^2 \rf^ \frac{1}{2} \leq ch 
\Vert u \Vert_{2,\infty,\bar \Omega},
\end{equation}
\begin{equation}
\label{T7E3}
\| e\|_{L^2} \leq ch \Vert u \Vert_{2,\infty,\bar \Omega},
\end{equation}
where $ e_{ij}=u(x_i,y_j)-u_{ij} $ for 
$(i,j) \in \{1,...,M\}\times \{1,...,N\}$ and vanishes elsewhere.
\end{theorem} 
\subsection{The First Correction} We proceed as in the third section, we begin 
by finding an expansion of the error. Begining by $\int_{K_{ij}}pudxdy$\\
We have
\begin{equation}
\label{Int_pu}
\int_{K_{ij}}pudxdy=h_ik_jp_{ij}u(x_i,y_j)+k_j \frac{{h_i^{+}}^2-{h_i^{-}}^2}{2}\frac{\partial (pu)}{\partial x}(a_{ij})
+h_i \frac{{k_j^{+}}^2-{k_j^{-}}^2}{2}\frac{\partial (pu)}{\partial y}(a_{ij})+N_{ij}^1,
\end{equation} 
where $a_{ij}=(x_i,y_j)$ and
\begin{equation}
\label{NIJ1}
\vert N_{ij}^1\vert \leq c h_ik_jh^2 \Vert u \Vert_{2,\infty,\bar \Omega}.
\end{equation}
Combining equalities (\ref{EQ112}) and (\ref{Int_pu}) yields
\begin{equation}
\label{eq12}
\begin{array}{rl}
(&\Dhh u)_{ij}\dsp{+h_ik_j p_{ij}u(x_i,y_j)=
\int_{K_{ij}}fdxdy}\\
&\\
&-k_j\frac{h_{i+1}^{-}-h_i^{+}}{2}u_{2x}(x_{i+ \frac{1}{2}},y_j)+
k_j\frac{h_i^{-}-h_{i-1}^{+}}{2}u_{2x}(x_{i- \frac{1}{2}},y_j)
-h_i\frac{k_{j+1}^{-}-k_j^{+}}{2}u_{2y}(x_i,y_{j+ \frac{1}{2}})+
h_i\frac{k_j^{-}-k_{j-1}^{+}}{2}u_{2y}(x_i,y_{j- \frac{1}{2}})\\
&\\
&+h_i \frac{{k_j^{+}}^2}{2}u_{2x,y}(x_i,y_j)-
h_i \frac{{k_j^{-}}^2}{2}u_{2x,y}(x_i,y_{j-1})
+k_j \frac{{h_i^{+}}^2}{2}u_{x,2y}(x_i,y_j)-
k_j \frac{{h_i^{-}}^2}{2}u_{x,2y}(x_{i-1},y_j) \\
&\\
&-k_j\frac{{h_i^{+}}^2-{h_i^{-}}^2}{2}\lf p_x(x_i,y_j)u(x_i,y_j)
+p_{ij}u_x(x_i,y_j)\rf 
-h_i\frac{{k_j^{+}}^2-{k_j^{-}}^2}{2}\lf p_y(x_i,y_j)u(x_i,y_j)
+p_{ij}u_y(x_i,y_j)\rf \\
&\\
&-k_j(R_{\IIj}^2-R_{\Iij }^2)         
+\frac{{k_j^{+}}^2-{k_j^{-}}^2}{2}(S_{ij}^1-S_{ij}^2) 
-h_i(R_{\iJJ}^2+R_{\iJj }^2)         
+\frac{{h_i^{+}}^2-{h_i^{-}}^2}{2}(L_{ij}^1-L_{ij}^2)\\
&\\
&+T_{ij}+H_{ij}+U_{ij}^1+U_{ij}^2-N_{ij}^1.
\end{array}
\end{equation}
We look, now, for an approximation to the second derivative $u_{2x}$ of $u$ by
using the same matrix that used to compute the basic solution $u^h$.
Remarking that $u_{2x}$ is the solution of the following equation
\begin{equation}
\label{Eq_DV}
-\Delta v+pv=f_{2x}-p_{2x}u-2p_xu_x
\end{equation} 
with the boundary conditions
\begin{equation}
\label{V_BORD} 
\left\lbrace
\begin{array}{l}
v|_ {\Gamma_1}=v|_ {\Gamma_3}=0,\\
v|_ {\Gamma_2}=-f|_ {\Gamma_2},\\
v|_ {\Gamma_4}=-f|_ {\Gamma_4}.
\end {array}
\right.
\end{equation}
Hence an approximation $v^h=(v_{ij})$ to $u_{2x}$ can be defined as
\begin{equation}
\label{V_INI} 
\left\lbrace
\begin{array}{l}
v_{i0}=v_{i,N+1}=0,\\
v_{0j}=-f(0,y_j),\\
v_{M+1,j}=-f(1,y_j)
\end {array}
\right.
\end{equation}
and for $(i,j) \in \{1,...,M\} \times \{1,...,N\}$, we have
\begin{eqnarray}
\label{SchV}
(\Dh v^h)_{ij}&+&h_ik_jp_{ij}v_{ij}=\int_{K_{ij}}f_{2x}dxdy\nonumber \\
&-&h_ik_jp_{2x}(x_i,y_j)u_{ij}-
2h_i^{+}k_jp_{x}(x_i,y_j)\partial_1^x u_{ij}
-2h_i^{-}k_jp_{x}(x_i,y_j)\partial_1^x u_{i-1,j}
\end{eqnarray}
where $u_{ij}$ are the components of the finite volume solution $u^h$ defined 
by (\ref{BorM})-(\ref{Sch_VFM}).\\
To analyse the convergence of the finite volume approximation $v^h$, let
$v=u_{2x}$ and using equality (\ref{VF_PM})
to get
\begin{eqnarray}
(\Dhh u_{2x})_{ij}&+&h_ik_jv(x_i,y_j)=
\int_{K_{ij}}(f-p_{2x}u-2p_xu_x)dxdy\nonumber \\
&+&k_j\lf R_{\IIj}(u_{2x})-R_{\Iij }(u_{2x})\rf+
h_i\lf R_{\iJJ}(u_{2x})-R_{\iJj }(u_{2x})\rf\nonumber \\
&+&S_{\IIj}(u_{2x})-S_{\Iij }(u_{2x})
+S_{\iJJ}(u_{2x})-S_{\iJj }(u_{2x})-N_{ij}(u_{2x}),
\end{eqnarray}
where $R_{\IIj}(u_{2x})$, $R_{\iJJ}(u_{2x})$, $S_{\IIj}(u_{2x})$, 
$S_{\iJJ}(u_{2x})$, $N_{ij}(u_{2x})$ are the same previous expansions by 
substituting each $u$ by $u_{2x}$.\\
This implies that
\begin{eqnarray}
\label{SchV1}
(\Dhh u_{2x})_{ij}&+&h_ik_jv(x_i,y_j)=\int_{K_{ij}}fdxdy\nonumber \\
&-&h_ik_jp_{2x}(x_i,y_j)u(x_i,y_j)-2k_jh_i^{+}p_x(x_i,y_j)u_x(x_i,y_j)\nonumber \\
&-&2k_jh_i^{-}p_x(x_i,y_j)u_x(x_{i-1},y_j)+O(h_ik_jh)+O(h_i^{-}k_j\Hi )\nonumber \\
&+&k_jR_{\IIj}(u_{2x})-k_jR_{\Iij }(u_{2x})+h_iR_{\iJJ}(u_{2x})-h_iR_{\iJj }(u_{2x})\nonumber \\
&+&S_{\IIj}(u_{2x})-S_{\Iij }(u_{2x})+S_{\iJJ}(u_{2x})-S_{\iJj }(u_{2x})-N_{ij}(u_{2x})
\end{eqnarray}
Let $r^h=(r_{ij})_{i,j}=(v(x_i,y_j)-v_{ij})_{i,j}$ be the error in the 
approxomation (\ref{SchV}). Substracting (\ref{SchV}) from (\ref{SchV1}), we
get for $(i,j) \in \{1,...,M\} \times \{1,...,N\}$
\begin{eqnarray}
\label{SchR}
(\Dh r)_{ij}+h_ik_jp_{ij}r_{ij} 
=\alpha_{\IIj}-\alpha_{\Iij }
+\alpha_{\iJJ}-\alpha_{\iJj }
+\delta_{ij}
\end{eqnarray}
where
\begin{eqnarray}
\alpha_{\IIj}=k_jR_{\IIj}(u_{2x})+S_{\IIj}(u_{2x}),\\
\alpha_{\iJJ}=h_iR_{\iJJ}(u_{2x})+S_{\iJJ}(u_{2x}),
\end{eqnarray}
\begin{equation}
\left\{\begin{array}{rcl}
\delta_{ij}&=&-h_ik_jp_{2x}(x_i,y_j)(u(x_i,y_j)-u_{ij})-2k_jh_i^{+}p_x(x_i,y_j)(u_x(x_i,y_j)-\partial_1^xu_{ij}) \\
&&\\
&&-2k_jh_i^{-}p_x(x_i,y_j)(u_x(x_{i-1},y_j)-\partial_1^xu_{i-1,j}) +O(h_ik_jh)+O(h_i^{-}k_j\Hi ).
\end{array}
\right.
\end{equation} 
Multiplying both sides of (\ref{SchR}) by $r_{ij}$ and summing over 
$(i,j) \in \{1,...,M\} \times \{1,...,N\}$, to get
\begin{equation}
\label{SR_ij}
\|r^h\|_{1,\tau}^2+\sMN h_ik_jp_{ij}r_{ij}^2
=-\SMN \alpha_{\IIj}(r_{i+1,j}-r_{ij})-
\SMN \alpha_{\IIj}(r_{i,j+1}-r_{ij})
+\sMN \delta_{ij} r_{ij}.
\end{equation}
Using the tricks those used to bound the error in the first correction
(subsection 3.3) to obtain the following optimal approximation to $u_{2x}$
\begin{lemma}
\label{Lemma9}
Let $u$ be the solution of the equation (\ref{pb2dM}). If 
$u \in C^4(\bar \Omega)$ and $p \in C^2(\bar \Omega)$. Then the approximation 
$v^h$ defined by the boundary conditions (\ref{V_INI}) and the equation 
(\ref{SchV}) satifies the following estimates
\begin{eqnarray}
\lf\SMN h_ik_j(u_{2x}(x_{\II },y_j)-v_{ij})^2\rf^ \frac{1}{2} 
\leq ch \Vert u \Vert_{4,\infty,\bar \Omega}
\end{eqnarray}
\begin{eqnarray}
\lf\SMN \HI k_j(u_{2x}(x_{\II },y_j)-v_{ij})^2\rf^ \frac{1}{2} 
\leq ch \Vert u \Vert_{4,\infty,\bar \Omega}
\end{eqnarray}
\begin{eqnarray}
\lf\SMN \KJ h_i(u_{2x}(x_{\II },y_j)-v_{ij})^2\rf^ \frac{1}{2} 
\leq ch \Vert u \Vert_{4,\infty,\bar \Omega}
\end{eqnarray}
\begin{eqnarray}
\lf\SMN h_i\KJ (u_{2x,y}(x_i,y_j)-\frac{v_{i,j+1}-v_{ij}}{\KJ })^2\rf^ \frac{1}{2} 
\leq ch \Vert u \Vert_{4,\infty,\bar \Omega}
\end{eqnarray}
\begin{eqnarray}
\lf\SMN \HI k_j(u_{3x}(x_i,y_j)-\frac{v_{i+1,j}-v_{ij}}{\HI })^2\rf^ \frac{1}{2} 
\leq ch \Vert u \Vert_{4,\infty,\bar \Omega}
\end{eqnarray}
\end{lemma} 
After having acheived an optimal approximation for $u_{2x}$ , we have to find 
an improvement of the basic solution $u^h$, i.e. correction of order $O(h^2)$. 
Looking again at the equality (\ref{eq12}) and rewrite
\begin{eqnarray}
(\Dhh u)_{ij}&+&h_ik_jp_{ij}u(x_i,y_j)=
\int_{K_{ij}}fdxdy
-k_j\frac{h_{i+1}^{-}-h_i^{+}}{2}u_{2x}(x_{i+ \frac{1}{2}},y_j)+
k_j\frac{h_i^{-}-h_{i-1}^{+}}{2}u_{2x}(x_{i- \frac{1}{2}},y_j)\nonumber\\
&&-h_i\frac{k_{j+1}^{-}-k_j^{+}}{2}u_{2y}(x_i,y_{j+ \frac{1}{2}})+
h_i\frac{k_j^{-}-k_{j-1}^{+}}{2}u_{2y}(x_i,y_{j- \frac{1}{2}})\nonumber \\ 
&&+h_i \frac{{k_j^{+}}^2}{2}u_{2x,y}(x_i,y_j)-
h_i \frac{{k_j^{-}}^2}{2}u_{2x,y}(x_i,y_{j-1})+k_j \frac{{h_i^{+}}^2}{2}u_{x,2y}(x_i,y_j)-
k_j \frac{{h_i^{-}}^2}{2}u_{x,2y}(x_{i-1},y_j)\nonumber \\
&&-k_j \frac{{h_i^{+}}^2-{h_i^{-}}^2}{2}p_x(x_i,y_j)u(x_i,y_j)
-k_j \frac{{h_i^{+}}^2}{2}p_{ij}u_x(x_i,y_j)\nonumber \\
&&+k_j \frac{{h_i^{-}}^2}{2}p_{ij}u_x(x_{i-1},y_j)-
h_i \frac{{k_j^{+}}^2-{k_j^{-}}^2}{2}p_y(x_i,y_j)u(x_i,y_j)\nonumber \\
&&-h_i \frac{{k_j^{+}}^2}{2}p_{ij}u_y(x_i,y_j)+
h_i \frac{{k_j^{-}}^2}{2}p_{ij}u_y(x_i,y_{j-1})\nonumber \\
&&-k_j(R_{\IIj}^2-R_{\Iij }^2)         
+\frac{{k_j^{+}}^2-{k_j^{-}}^2}{2}(S_{ij}^1-S_{ij}^2)+T_{ij}\nonumber \\
&&-h_i(R_{\iJJ}^2-R_{\iJj }^2)         
+\frac{{h_i^{+}}^2-{h_i^{-}}^2}{2}(L_{ij}^1-L_{ij}^2)+H_{ij}\nonumber \\
&&+U_{ij}^1+U_{ij}^2-N_{ij}^1+A_{ij}+B_{ij}
\end{eqnarray} 
where
\begin{equation}
\label{EstAIJ}
|A_{ij}|\leq ck_j {h_i^{-}}^2 \Hi  \Vert u \Vert_{2,\infty,\bar \Omega},
\end{equation}
\begin{equation}
\label{EstBIJ}
|B_{ij}|\leq ch_i {k_j^{-}}^2 k_{\Jj } \Vert u \Vert_{2,\infty,\bar \Omega}.
\end{equation} 
To simplify the expressions, let
\begin{equation}
\label{bIJ}
b_{ij}=-k_j\frac{h_{i+1}^{-}-h_i^{+}}{2}v_{ij}
\end{equation}
\begin{equation}
\label{cIJ}
c_{ij}=-h_i\frac{h_{i+1}^{-}-h_i^{+}}{2}\lf -f(x_i,y_{\JJ })+
p(x_i,y_{\JJ })u_{ij}-v_{ij}\rf
\end{equation}
\begin{eqnarray}
d_{ij}&&=h_i \frac{{k_j^{+}}^2}{2}\partial_1^y v_{ij}-h_i \frac{{k_j^{-}}^2}{2}\partial_1^y v_{i,j-1}+
k_j \frac{{h_i^{+}}^2}{2}\lf-f_x(x_i,y_j)+p_x(x_i,y_j)u_{ij}+p_{ij}\partial_1^x u_{ij}-\partial_1^x v_{ij}\rf
\nonumber \\
&&-k_j \frac{{h_i^{-}}^2}{2}\lf -f_x(x_{i-1},y_j)+p_x(x_{i-1},y_j)u_{i-1,j}+p_{i-1,j}
\partial_1^x u_{i-1,j}-\partial_1^x v_{i-1,j}\rf \nonumber \\
&&-k_j \frac{{h_i^{+}}^2-{h_i^{-}}^2}{2}p_x(x_i,y_j)u_{ij}
-k_j \frac{{h_i^{+}}^2}{2}p_{ij}\partial_1^xu_{ij}+
k_j \frac{{h_i^{-}}^2}{2}p_{ij}\partial_1^xu_{i-1,j}\nonumber \\
&&-h_i \frac{{k_j^{+}}^2-{k_j^{-}}^2}{2}p_y(x_i,y_j)u_{ij}
-h_i \frac{{k_j^{+}}^2}{2}p_{ij}\partial_1^yu_{ij}+
h_i \frac{{k_j^{-}}^2}{2}p_{ij}\partial_1^yu_{i,j-1}
\end{eqnarray} 
Now, the first correction $u^h_1=(u_{ij}^1)$ can be defined as follows
\begin{equation}
\label{FirstBor}
u^1_{0j}=u^1_{M+1,j}=u^1_{i0}=u^1_{i,N+1}=0,
\end{equation}
and for $(i,j) \in \{1,...,M\} \times \{1,...,N\}$, we have
\begin{eqnarray}
(\Dh u^h_1)_{ij}+h_ik_jp_{ij}u_{ij}=\int_{K_{ij}}fdxdy
+ b_{ij}-b_{i-1,j}+c_{ij}-c_{i,j-1}+d_{ij}.
\end{eqnarray}
Using inequalities (\ref{Est_TIJp}), (\ref{Est_SIJ1}), (\ref{Est_SIJp}),
(\ref{Est_UIJ1}),
(\ref{Est_RIJ2}),(\ref{Est_RJ2})-(\ref{Est_UIJ2}), (\ref{NIJ1}),
(\ref{EstAIJ}) and (\ref{EstBIJ})
combined with the estimates obtained in lemma \ref{Lemma9} yields the following
$O(h)$ improvement
\begin{theorem}
If the solution $u$ of the equation (\ref{pb2dM}) belonging to 
$C^4(\bar \Omega)$ and $p \in C^2(\bar \Omega)$. Let $u^h$ be the basic finite 
volume solution of boundary condition (\ref{BorM}) and discrete
equation (\ref{Sch_VFM}). Then the finite volume approximations 
$u^h_1=(u_{ij}^1)$ defined by the boundary condition (\ref{FirstBor}) and the 
discrete equation (\ref{FirstBor}), satisfies the following estimates
\begin{equation}
\|e^1\|_{1,\tau} \leq ch^2 \Vert u \Vert_{4,\infty,\bar \Omega},
\end{equation}
\begin{equation}
\lf\SMN k_j\HI (e_{ij}^1)^2\rf^ \frac{1}{2} \leq ch^2 
\Vert u \Vert_{4,\infty,\bar \Omega},
\end{equation}
\begin{equation}
\|e^1\|_{L^2} \leq ch^2 \Vert u \Vert_{4,\infty,\bar \Omega},
\end{equation}
where $ e_{ij}^1=u(x_i,y_j)-u_{ij}^1 $ for $(i,j) \in \{1,...,M\}\times \{1,...,N\}$ and
vanishes elsewhere.
\end{theorem}


\newpage
\section{Numerical Tests}
\subsection{In one Dimensional Space}
In this subsection, we give two numerical tests justifying our theoretical
results in one dimensional case. We mean by uniform mesh that so-called 
modified finite volume scheme and satisfying $\HI =h $, $h_{i+1}^{-}=h_i^{+}$ 
(see remark 2.5 in \cite{Eymard} ) , and by cell-centered mesh 
that satisfying 
$h_i= \left\lbrace 
\begin{array}{l}
h, i  \mbox{ is  even}, \\
\frac{h}{2}, i  \mbox{ is  odd}
\end{array} 
\right.
$
and $h_i^{-}=h_i^{+}.$  Note that, the first correction computed in case of 
uniform mesh is the second one, because  $h_{i+1}^{-}-h_i^{+}=0$ .\\
 To show the convergence orders of the first correction and the basic finite
volume solution, we compute the ratio
$$
ratio=\frac{\log (e(h))-\log (e(h_0))}{\log (h)-\log (h_0)}.
$$
where $h_0$ is the initial value of $h$ in each numerical test and $e(h)$ is 
the error corresponding to $h$. In the uniform mesh, we use the rule 
$$
ratio=-\frac{\log (e(1/2^{k+1}))-\log (e(1/2^k))}{\log 2}.
$$
\subsubsection{{\bf{First Test}}} We consider the homogeneous equation
${\bf{(I)}}: -u_{xx}=f$ 
where $ u(x)=sin(\pi x) $ and $f(x)=\pi^2sin(\pi x)$.
\begin{center}
{\bf{TABLE 1.}} The convergence orders of the first correction and the basic solution in
$L^2$-norm in uniform mesh.
\end{center}
\begin{center}
\begin{tabular}{|c||c|c||c|c||}
\hline
&
\multicolumn{2}{|c||}{}
&
\multicolumn{2}{|c||}{}\\
$h$
&
\multicolumn{2}{|c||}{\mbox{correction}}
&
\multicolumn{2}{|c||}{\mbox{basic solution}}\\
\cline{2-5}

&order&error /$h^4$&order&error /$h^4$\\
\hline
\multicolumn{1}{|c||}{1/32}
&-&0.0359&-&0.0298e+04 \\
\hline
1/64&4.0004&0.0359&2.0003&0.1191e+04 \\
\hline
1/128&4.0000&0.0359&2.0001&0.4764e+04\\
\hline
1/256&4.0032&0.0390&2.0000&1.9057e+04\\
\hline
\end{tabular}
\end{center}
\begin{center}
{\bf{TABLE 2. }} The convergence orders of the first correction and the basic solution in
$H^1_0$-norm in uniform mesh.
\end{center}
\begin{center}
\begin{tabular}{|c||c|c||c|c||}
\hline
&
\multicolumn{2}{|c||}{}
&
\multicolumn{2}{|c||}{}\\
$h$
&
\multicolumn{2}{|c||}{\mbox{correction}}
&
\multicolumn{2}{|c||}{\mbox{basic solution}}\\
\cline{2-5}

&order&error /$h^4$&order&error /$h^4$\\
\hline
\multicolumn{1}{|c||}{1/32}
&-&0.1127&-&0.0935e+04 \\
\hline
1/64&3.9999&0.1127&1.9999&0.3742e+04 \\
\hline
1/128&3.9999&0.1129&2.0000&1.4967e+04\\
\hline
1/256&4.0032&0.1224&2.0000&5.9869e+04\\
\hline
\end{tabular}
\end{center}

\begin{center}
{\bf{TABLE 3.}} The convergence orders of the first correction and the basic solution in
$L^2$-norm in cell-centered mesh.
\end{center}
\begin{center}
\begin{tabular}{|c||c|c||c|c||}
\hline
&
\multicolumn{2}{|c||}{}
&
\multicolumn{2}{|c||}{}\\
$h$
&
\multicolumn{2}{|c||}{\mbox{correction}}
&
\multicolumn{2}{|c||}{\mbox{basic solution}}\\
\cline{2-5}

&order&error /$h^2$&order&error /$h^2$\\
\hline
\multicolumn{1}{|c||}{4/149}
&-&0.2474&-&0.3971\\
\hline
4/599&1.9903&0.2507&1.9986&0.3978\\
\hline
4/2999&1.9943&0.2516&1.9991&0.3981\\
\hline
4/14999&1.9961&0.2473&1.9994&0.3968\\
\hline
\end{tabular}
\end{center}
\begin{center}
{\bf{TABLE 4. }} The convergence orders of the first correction and the basic solution in
$H^1_0$-norm in cell-centered mesh.
\end{center}
\begin{center}
\begin{tabular}{|c||c|c||c|c||}
\hline
&
\multicolumn{2}{|c||}{}
&
\multicolumn{2}{|c||}{}\\
$h$
&
\multicolumn{2}{|c||}{\mbox{correction}}
&
\multicolumn{2}{|c||}{\mbox{basic solution}}\\
\cline{2-5}

&order&error /$h^2$&order&error /$h^2$\\
\hline
\multicolumn{1}{|c||}{4/149}
&-&0.5856&-&0.0281e+03\\
\hline
4/599&1.9929&0.5914&1.0000&0.1131e+03\\
\hline
4/2999&1.9958&0.5930&1.0000&0.5664e+03\\
\hline
4/14999&1.9971&0.5828&1.0000&2.8329e+03\\
\hline
\end{tabular}
\end{center}
\begin{center}
{\bf{TABLE 5. }} Comparaison between the accuracy of the first correction that
uses first variant and the one using second variant in $ H^1_0$ and $L^2$-norms 
in cell-centered mesh.
\end{center}
\begin{center}
\begin{tabular}{|c||c|c||c|c||}
\hline
&
\multicolumn{2}{|c||}{}
&
\multicolumn{2}{|c||}{}\\
$h$
&
\multicolumn{2}{|c||}{\mbox{first variant}}
&
\multicolumn{2}{|c||}{\mbox{second variant}}\\
\cline{2-5}

&$L^2$-norm&$H^1_0$-norm&$L^2$-norm&$H^1_0$-norm\\
\hline
\multicolumn{1}{|c||}{4/149}
&1.7827e-04&4.2204e-04&7.9569e-05&6.6012e-04\\
\hline
4/599&1.1180e-05&2.6374e-05&5.4358e-06&4.1591e-05\\
\hline
4/2999&4.4758e-07&1.0549e-06&2.2254e-07&1.6670e-06\\
\hline
4/14999&1.7591e-08&4.1450e-08&8.6276e-09&6.6477e-08\\
\hline
\end{tabular}
\end{center}
\subsubsection{{\bf{Second Test}}} In case of cell-centered mest, we saw that
the convergence of the first correction is the same one of the basic solution in
$L^2$-norm for the model ({\bf{I}}). We present here an example of the mesh where the convergence order of the first
correction improves really that one of the basic solution in both $H^1_0$ and
$L^2$-norms for the model ({\bf{I}}).\\
We consider $h_\II=h$ and $x_\II=\frac{2x_i+x_{i+1}}{3}$ for all $i=1,..., N-1$.
\begin{center}
{\bf{TABLE 6. }}The convergence orders of the first correction and the basic solution 
in $L^2$-norm .
\end{center}
\begin{center}
\begin{tabular}{|c||c|c||c|c||}
\hline
&
\multicolumn{2}{|c||}{}
&
\multicolumn{2}{|c||}{}\\
$h$
&
\multicolumn{2}{|c||}{\mbox{correction}}
&
\multicolumn{2}{|c||}{\mbox{basic solution}}\\
\cline{2-5}

&order&error /$h^2$&order&error /$h^2$\\
\hline
\multicolumn{1}{|c||}{$1/2^8$}
&-&0.3701&-&20.2232\\
\hline
$1/2^9$&1.9668&0.3788&1.0000&40.4472\\
\hline
$1/2^{10}$&1.9633&0.3832&1.0000&80.8948\\
\hline
$1/2^{11}$&1.9915&0.3855&1.0000&161.7899\\
\hline
$1/2^{12}$&2.0073&0.3835&1.0000&323.5799\\
\hline
$1/2^{13}$&1.9864&0.3871&1.0000&647.1598\\
\hline
\end{tabular}
\end{center}
\begin{center}
{\bf{TABLE 7. }} The convergence orders of the first correction and the basic solution in
$H^1_0$-norm.
\end{center}
\begin{center}
\begin{tabular}{|c||c|c||c|c||}
\hline
&
\multicolumn{2}{|c||}{}
&
\multicolumn{2}{|c||}{}\\
$h$
&
\multicolumn{2}{|c||}{\mbox{correction}}
&
\multicolumn{2}{|c||}{\mbox{basic solution}}\\
\cline{2-5}

&order&error /$h^2$&order&error /$h^2$\\
\hline
\multicolumn{1}{|c||}{$1/2^8$}
&-&1.4588&-&129.5822\\
\hline
$1/2^9$&2.1180&1.3443&0.9999&259.1881\\
\hline
$1/2^{10}$&2.0675&1.2828&1.0000&518.3882\\
\hline
$1/2^{11}$&2.0363&1.2509&1.0000&1.0368e+03\\
\hline
$1/2^{12}$&2.0301&1.2251&1.0000&2.0736e+03\\
\hline
$1/2^{13}$&1.9985&1.2263&1.0000&4.1471e+03\\
\hline
\end{tabular}
\end{center}
\subsubsection{{\bf{Third Test}}}We consider the homogeneous equation 
$
{\bf{(II)}}: -u_{xx}+u_x+u=f
$
where $u(x)=sin(\pi x) $ and $f(x)=(\pi^2+1)\sin(\pi x)+\pi \cos (\pi x).$
\begin{center}
{\bf{TABLE 8. }}The convergence orders of the first correction and the basic solution 
in $L^2$-norm in cell-centered mesh.
\end{center}
\begin{center}
\begin{tabular}{|c||c|c||c|c||}
\hline
&
\multicolumn{2}{|c||}{}
&
\multicolumn{2}{|c||}{}\\
$h$
&
\multicolumn{2}{|c||}{\mbox{correction}}
&
\multicolumn{2}{|c||}{\mbox{basic solution}}\\
\cline{2-5}

&order&error /$h^2$&order&error /$h^2$\\
\hline
\multicolumn{1}{|c||}{2/75}
&-&0.1415&-&0.0099e+03\\
\hline
1/150&2.0074&0.1401&0.9796&0.0407e+03\\
\hline
1/750&2.0043&0.1398&0.9880&0.2051e+03\\
\hline
1/3750&2.0020&0.1402&0.9918&1.0270e+03\\
\hline
\end{tabular}
\end{center}
\begin{center}
{\bf{TABLE 9. }} The convergence orders of the first correction and the basic solution in
$H^1_0$-norm in cell-centered mesh.
\end{center}
\begin{center}
\begin{tabular}{|c||c|c||c|c||}
\hline
&
\multicolumn{2}{|c||}{}
&
\multicolumn{2}{|c||}{}\\
$h$
&
\multicolumn{2}{|c||}{\mbox{correction}}
&
\multicolumn{2}{|c||}{\mbox{basic solution}}\\
\cline{2-5}

&order&error /$h^2$&order&error /$h^2$\\
\hline
\multicolumn{1}{|c||}{2/75}
&-&0.4441&-&0.0377e+03\\
\hline
1/150&1.9936&0.4480&0.9887&0.1530e+03\\
\hline
1/750&1.9963&0.4491&0.9934&0.7681e+03\\
\hline
1/3750&1.9970&0.4502&0.9955&3.8438e+03\\
\hline
\end{tabular}
\end{center}
\subsection{{\bf{In Two Dimensional Space }}} In this subsection, we present two
tests justifying our results of two dimensional space. The ratios are computed
by using the first formula of ratio in subsection 5.1 .   
\subsubsection{{\bf{First Test }}}
We consider here $u=xy(1-x)(1-y)$, then $u$
 is the solution of (\ref{pb2d}), where $f=2y(1-y)+2x(1-x)$. The mesh 
considered here is such that $h_i= \left\lbrace 
\begin{array}{l}
h, i  \mbox{ is  even}, \\
\frac{h}{2}, i  \mbox{ is  odd}
\end{array} 
\right.
$
and $k_j=3h/2$, with $(x_i,y_j)$ is in the center of $K_{ij}$.
\begin{center}
{\bf{TABLE 10.}} Convergence orders of the first correction and basic solution in 
$H^1_0$-norm.\\ 
\end{center}
\begin{center}
\begin{tabular}{|c||c|c||c|c||}
\hline
&
\multicolumn{2}{|c||}{}
&
\multicolumn{2}{|c||}{}\\
$h$
&
\multicolumn{2}{|c||}{\mbox{correction}}
&
\multicolumn{2}{|c||}{\mbox{basic solution}}\\
\cline{2-5}

&order&error /$h^2$&order&error /$h^2$\\
\hline
\multicolumn{1}{|c||}{0.02500}
&-&0.198028&-&2.51479\\
\hline
0.01667&1.99847&0.198151&1.02547&3.73344\\
\hline
0.01250&1.99906&0.198157&1.02262&4.95135\\
\hline
0.01000&1.99983&0.198139&1.02068&6.16896\\
\hline
0.00500&1.99992&0.198053&1.01593&12.2557\\
\hline
0.00333&2.00005&0.198007&1.01385&18.3418\\
\hline
\end{tabular}
\end{center}
\begin{center}
{\bf{TABLE 11. }}Convergence orders of the first correction and the basic 
solution in 
$L^2$-norm.
\end{center}
\begin{center}
\begin{tabular}{|c||c|c||c|c||}
\hline
&
\multicolumn{2}{|c||}{}
&
\multicolumn{2}{|c||}{}\\
$h$
&
\multicolumn{2}{|c||}{\mbox{correction}}
&
\multicolumn{2}{|c||}{\mbox{basic solution}}\\
\cline{2-5}

&order&error /$h^2$&order&error /$h^2$\\
\hline
\multicolumn{1}{|c||}{0.02500}
&-&0.0439562&-&0.07314\\
\hline
0.01667&1.94335&0.0445737&1.99851&0.07318\\
\hline
0.01250&1.95651&0.0448951&1.99882&0.07320\\
\hline
0.01000&1.96233&0.045092&1.99900&0.07321\\
\hline
0.00500&1.97302&0.0454949&1.99935&0.07322\\
\hline
0.00333&1.97708&0.0456318&1.99947&0.07322\\
\hline
\end{tabular}
\end{center}


\subsubsection{{\bf{Second Test}}}
In this test we choose $u(x,y)=\sin(\pi x) \sin(\pi y) $, then $u$ is the 
solution of the problem (\ref{pb2d}), where 
$f(x,y)=2\pi^2\sin(\pi x) \sin(\pi y)$. The mesh considered here is such that 
$(x_i,y_j)$ is in the center of $K_{ij}$, 
$h_i= \left\lbrace 
\begin{array}{l}
h, i  \mbox{ is  even}, \\
\frac{h}{2}, i  \mbox{ is  odd}
\end{array} 
\right.
$
and $k_j= \left\lbrace 
\begin{array}{l}
h, j  \mbox{ is  even}, \\
\frac{h}{2}, j  \mbox{ is  odd}
\end{array} 
\right.
$. 




\begin{center}
{\bf{TABLE 12. }} The convergence orders of the first correction and the basic solution in
$H^1_0$-norm.
\end{center}
\begin{center}
\begin{tabular}{|c||c|c||c|c||}
\hline
&
\multicolumn{2}{|c||}{}
&
\multicolumn{2}{|c||}{}\\
$h$
&
\multicolumn{2}{|c||}{\mbox{correction}}
&
\multicolumn{2}{|c||}{\mbox{basic solution}}\\
\cline{2-5}

&order&error /$h^2$&order&error /$h^2$\\
\hline
\multicolumn{1}{|c||}{0.033333}
&-&2.80465&-&34.946\\
\hline
0.016667&2.21288&2.4199&1.0016&69.815\\
\hline
0.008333&2.17818&2.1908&1.0010&139.59\\
\hline
0.004166&2.14753&2.06364&1.0007&279.16\\
\hline
\end{tabular}
\end{center}
\begin{center}
{\bf{TABLE 13.}} The convergence orders of the first correction and the basic solution in
$L^2$-norm.
\end{center}
\begin{center}
\begin{tabular}{|c||c|c||c|c||}
\hline
&
\multicolumn{2}{|c||}{}
&
\multicolumn{2}{|c||}{}\\
$h$
&
\multicolumn{2}{|c||}{\mbox{correction}}
&
\multicolumn{2}{|c||}{\mbox{basic solution}}\\
\cline{2-5}

&order&error /$h^2$&order&error /$h^2$\\
\hline
\multicolumn{1}{|c||}{0.033333}
&-&0.14632&-&0.7812\\
\hline
0.016667&2.15317&0.131583&2.00058&0.7809\\
\hline
0.008333&2.05991&0.134658&2.00036&0.7802\\
\hline
0.004166&2.02539&0.138795&2.00025&0.7802\\
\hline
\end{tabular}
\end{center}
\subsection{Some Comments about the Numerical Results. }
\begin{enumerate}
\item In {\bf{Table 1}} and {\bf{Table 2}}, numerical results show that on uniform mesh and for
the model {\bf{(I)}}, we can gain an $O(h^2)$-
improvement in both $H^1_0$-norm and $L^2$-norm by the first correction.
\item In {\bf{Table 3}}, numerical results show that for the model {\bf{(I)}}, 
we do not have an improvement in
$L^2$-norm in the first correction when the mesh cell-centered. Furthermore, 
the coefficents of the error in the first correction are 
better than of those of the basic solution. To improve the order in $L^2$-norm, 
we compute the second correction.
\item In {\bf{Table 4}}, numerical results show that for the model {\bf{(I)}}, we gain an $O(h)$-improvement
by the first correction in $H^1_0$-norm.
\item In {\bf{Table 5}}, numerical results show that the accuracy of the error 
in the correction defined by second variant (see subsection 2.3.1) is batter
than that of (\ref{1deq22}) in $L^2$-norm and contrary in $H^1_0$-norm.
\item In  {\bf{Table 6}} and {\bf{Table 7}}, numerical results show that the
convergence of the first correction improves that of the basic solution in both 
$H^1_0$ and $L^2$ norms for {\bf{(I)}}.\\
This implies that, on arbitrarily admissible mesh, the first correction improves
the basic solution in $L^2$ and $H^1_0$ norms. 
\item In  {\bf{Table 8}} and {\bf{Table 9}}, numerical results show that for the
model {\bf{(II)}}, we gain an $O(h)$-improvement in both $L^2$-norm and
$H^1_0$-norm by the first correction for cell-centered mesh.
\item In {\bf{Table 10}} and {\bf{Table 12}}, numerical results show that we gain
an $O(h)$-improvement in $H^1_0$-norm by the first correction.
\item In {\bf{Table 11}} and {\bf{Table 13}}, numerical results show that,
the convergence order of the first correction is the same one as of the basic
solution in $L^2$-norm, but the errors in the first correction are better than
of those of the basic solution.
\end{enumerate}   

\begin{remarque}. The idea used here is used by the authors to apply the defect correction technique in finite element method with non uniform mesh \cite{Atfeh}\end{remarque}
{\bf{Acknowledgment.}} The authors would like to thank Professor T.
Gallou\"{e}t who has attracted their attention to this worthy topic and his
sincere guidence and useful discussions


\begin{thebibliography}{12}
\bibitem{Atfeh}{\sc B. Atfeh and A. Bradji: }{Defect Correction in Finite Element Method with Non Uniform Mesh. }{\it In Preparation.}
\bibitem{Barrett}{\sc J. W. Barrett, G. Moore: }{Optimal Recovery in the Finite Element Method, Part 2 Defect Correction for Ordinary Differential Equations. }{\it IMA J.
Numer. Anal.},527-540, {\bf{8}}, 1988.
\bibitem{Bohmer}{\sc K. Bohmer and H. J. Stetter (eds):}{Defect Correction Theory and Applications. }{\it Springer-Verlag. Wien, New York}, 1984.
\bibitem{Butcher}{\sc J. C. Butcher, J. R. Cash, G. Moore and R. D. Russell }{Defect Correction for
two-Point Boundary Value Problems on Nonequidistant Mesh.}{\it Math. Comp. }, 629-648,{\bf{64}}, 1995.
\bibitem{Cai}{\sc Z. Cai, J. Douglas and M. Park: }{ Development and Analysis of Higher Order Finite Volume Methods Over Rectangles for Elliptic Problems. }{\it Advances in Computational Mathematics}, 3-33, {\bf{19}},2003.
\bibitem{Chibi}{\sc A. S.-Chibi : }{Defect Correction and Galerkin's Method for Second Order Elliptic Boundary Value Problems. }{\it Ph.D Thesis, Imperial college, London},1989.
\bibitem{Eymard}{\sc R. Eymard, T. Gallou\"{e}t and R. Herbin :}{ Finite Volume Methods. }{\it Handbook of Numerical Analysis. P. G. Ciarlet and J. L. Lions (eds.)},{\bf{vol. VII}}, 723-1020,2000.
\bibitem{Forsyth}{\sc P. A. Forsyth, Jr. and P. H. Sammon : }{Quadratic Convergence for Cell-Centered Grids. }{\it Applied Numerical Mathematics}, 
377-394, {\bf{4}}, 1988.
\bibitem{Fox}{\sc L. Fox : }{Some improvements in the Use of Relaxation Methods for the Solution of Ordinary and Partial Defferential Equations. }{\it Proc. Roy. Soc. Lon Ser. A}, 31-59, {\bf{190}},1947.
\bibitem{Gao}{\sc Jun-bin Gao, Yi-du Yang and T. M. Shih : }{The Defect Iteration of the Finite Element for Elliptic Boundary Value Problems and
Petrov-Galerkin Approximation. }{\it J. Computational Mathematics}, 152-164,
{\bf{16}}, 1998.
\bibitem{Martin}{\sc R. Martin and H. Guillard : }{A Second Order Defect
Correction Scheme for Unsteady Problems. }{\it Rapport de rechereche INRIA.}, {\bf{2447}},1994.
\bibitem{Moore}{\sc G. Moore : }{Defect Correction from a Galerkin View Points .}
{\it Num. Math }, 565-582,
{\bf{52}}, 1988.
\bibitem{Pereyra}{\sc V. Pereyra : }{Iterated Defrred Corrections for non Linear Operator Equations}{\it Num. Math.}, 316-323, {\bf{10}}, 1967.
\bibitem{Skeel}{\sc R. D. Skeel : }{A Theoretical  Framework for Proving
Accuracy Results for Deferred Corrections. }{\it Siam J. Num. Anal}, 171-196, {\bf{19}},1981. 
 

\end{thebibliography}
\end{document}